\documentclass[12pt]{article}
\usepackage{latexsym,amssymb,hyperref}
\newtheorem{Def1}{Definition}
\newtheorem{Lemma}[Def1]{Lemma}
\newtheorem{Prop}[Def1]{Proposition}
\newtheorem{Theorem}[Def1]{Theorem}
\newtheorem{Cor}[Def1]{Corollary}

\newtheorem{Conj}[Def1]{Conjecture}

\newtheorem{Ex1}[Def1]{Example}
\newenvironment{Def} {\begin{Def1} \begin{upshape}} {\end{upshape} \end{Def1}}
\newenvironment{Example} {\begin{Ex1} \begin{upshape}} {\end{upshape} \end{Ex1}}

\def\defaultskip{\medskip}
\newenvironment{Proof}[1][.]{\defaultskip\noindent \textbf{Proof}#1 }
  {\hspace*{0mm}\hfill $\Box$ \defaultskip}
\newenvironment{keywords}{\centerline{\bf\small
Keywords}\begin{quote}\small}{\par\end{quote}\vskip 1ex}
\def\defaultskip{\medskip}

\def\beq{\begin{equation}}
\def\eeq{\end{equation}}
\def\beqn{\begin{displaymath}}
\def\eeqn{\end{displaymath}}
\def\bqa{\begin{eqnarray}}
\def\eqa{\end{eqnarray}}
\def\bqan{\begin{eqnarray*}}
\def\eqan{\end{eqnarray*}}

\def\NNN{\mathbb N}
\def\RRR{\mathbb R}
\def\QQQ{\mathbb Q}

\def\BBB{\mathbb B}

\def\Expect{{\mathbf E}}
\def\Prob{{\mathbf P}}

\def\eps{\varepsilon}

\def\equa{\stackrel+=}

\def\eqm{\stackrel\times=}
\def\leqm{\stackrel\times\leq}
\def\geqm{\stackrel\times\geq}

\def\leqn{_{1:n}}

\def\ltinf{_{<\infty}}
\def\_norm{_\mathrm{norm}}

\def\for_all{\mbox{ for all }}
\def\such_that{\mbox{ such that }}
\def\wenn{\mbox{ if }}
\def\und{\mbox{ and }}
\def\fuer{\mbox{ for }}

\def\lb{{\log_2}}

\def\l{\ell}

\def\ph{\varphi}
\def\th{\vartheta}

\def\eins{1\hspace{-0.23em}{\rm I}}

\def\M{M}

\def\K{K\!}
\def\Kpre{K}
\def\Km{K\!m}

\begin{document}

\title{\normalsize\sc Technical Report \hfill IDSIA-04-06
\vskip 2mm\bf\Large\hrule height5pt \vskip 6mm
MDL Convergence Speed for Bernoulli Sequences%
\thanks{A shorter version of this paper \cite{Poland:04mdlspeed} appeared in ALT 2004.}
\vskip 6mm \hrule height2pt}

\author{{\bf Jan Poland} \ and \ {\bf Marcus Hutter}\\[3mm]
IDSIA, Galleria 2 CH-6928 Manno (Lugano), Switzerland\thanks{This
work was supported by SNF grant 2100-67712.02.}\\
\texttt{\{jan,marcus\}@idsia.ch} \hspace{11ex} http://www.idsia.ch}

\date{22 February 2006}

\maketitle

\begin{abstract}
The Minimum Description Length principle for online sequence
estimation/prediction in a proper learning setup is studied. If
the underlying model class is discrete, then the total expected
square loss is a particularly interesting performance measure: (a)
this quantity is finitely bounded, implying convergence with
probability one, and (b) it additionally specifies the convergence
speed. For MDL, in general one can only have loss bounds which are
finite but exponentially larger than those for Bayes mixtures. We
show that this is even the case if the model class contains only
Bernoulli distributions. We derive a new upper bound on the
prediction error for countable Bernoulli classes. This implies a
small bound (comparable to the one for Bayes mixtures) for certain
important model classes. We discuss the application to Machine
Learning tasks such as classification and hypothesis testing, and
generalization to countable classes of i.i.d.\ models.
\end{abstract}

\vspace{4ex}

\newpage
\section{Introduction}\label{secIntro}

``Bayes mixture", ``Solomonoff induction", ``marginalization", all
these terms refer to a central induction principle: Obtain a
predictive distribution by integrating the product of prior and
evidence over the model class. In many cases however, the Bayes
mixture is computationally infeasible, and even a sophisticated
approximation is expensive. The MDL or MAP (maximum a posteriori)
estimator is both a common approximation for the Bayes mixture and
interesting for its own sake: Use the model with the largest
product of prior and evidence. (In practice, the MDL estimator is
usually being approximated too, since only a local maximum is
determined.)

How good are the predictions by Bayes mixtures and MDL? This
question has attracted much attention. In many cases, an important
quality measure is the \emph{total} or cumulative \emph{expected
loss} of a predictor. In particular the square loss is often
considered. Assume that the outcome space is finite, and the model
class is continuously parameterized. Then for Bayes mixture
prediction, the cumulative expected square loss is usually small
but unbounded, growing with $\ln n$, where $n$ is the sample size
\cite{Clarke:90,Hutter:03optisp}. This corresponds to an
\emph{instantaneous} loss bound of $\frac{1}{n}$. For the MDL
predictor, the losses behave similarly
\cite{Rissanen:96,Barron:98} under appropriate conditions, in
particular with a specific prior. (Note that in order to do MDL
for continuous model classes, one needs to \emph{discretize} the
parameter space, see also \cite{Barron:91}.)

On the other hand, if the model class is discrete, then
Solomonoff's theorem \cite{Solomonoff:78,Hutter:01alpha} bounds
the cumulative expected square loss for the Bayes mixture
predictions finitely, namely by $\ln w_\mu^{-1}$, where $w_\mu$ is
the prior weight of the ``true" model $\mu$. The only necessary
assumption is that the true distribution $\mu$ is contained in the
model class, i.e.\ that we are dealing with \emph{proper
learning}. It has been demonstrated \cite{Gruenwald:04}, that for
both Bayes mixture and MDL, the proper learning assumption can be
essential: If it is violated, then learning may fail very badly.

For MDL predictions in the proper learning case, it has been shown
\cite{Poland:04mdl} that a bound of $w_\mu^{-1}$ holds. This bound
is exponentially larger than the Solomonoff bound, and it is sharp
in general. A finite bound on the total expected square loss is
particularly interesting:
\begin{enumerate}
\item It implies convergence of the predictive to the true probabilities
with probability one. In contrast, an instantaneous loss bound of
$\frac{1}{n}$ implies only convergence in probability.
\item Additionally, it gives a \emph{convergence speed}, in the
sense that errors of a certain magnitude cannot occur too often.
\end{enumerate}
So for both, Bayes mixtures and MDL, convergence with probability
one holds, while the convergence speed is exponentially worse for
MDL compared to the Bayes mixture. (We avoid the term
``convergence rate" here, since the order of convergence is
identical in both cases. It is e.g.\ $o(1/n)$ if we additionally
assume that the error is monotonically decreasing, which is not
necessarily true in general).

It is therefore natural to ask if there are model classes where
the cumulative loss of MDL is comparable to that of Bayes mixture
predictions. In the present work, we concentrate on the simplest
possible stochastic case, namely discrete Bernoulli classes. (Note
that then the MDL ``predictor" just becomes an estimator, in that
it estimates the true parameter and directly uses that for
prediction. Nevertheless, for consistency of terminology, we keep
the term predictor.) It might be surprising to discover that in
general the cumulative loss is still exponential. On the other
hand, we will give mild conditions on the prior guaranteeing a
small bound. Moreover, it is well-known that the instantaneous
square loss of the Maximum Likelihood estimator decays as
$\frac{1}{n}$ in the Bernoulli case. The same holds for MDL, as we
will see. (If convergence speed is measured in terms of
instantaneous losses, then much more general statements are
possible \cite{Li:99,Zhang:04}, this is briefly discussed in
Section \ref{secUB}.)

A particular motivation to consider discrete model classes arises
in Algorithmic Information Theory. From a computational point of
view, the largest relevant model class is the class of all
computable models on some fixed universal Turing machine,
precisely prefix machine \cite{Li:97}. Thus
each model corresponds to a program, and there are countably many
programs. Moreover, the models are stochastic, precisely they are
\emph{semimeasures} on strings (programs need not halt, otherwise
the models were even measures). Each model has a natural
description length, namely the length of the corresponding
program. If we agree that programs are binary strings, then a
prior is defined by two to the negative description length. By the
Kraft inequality, the priors sum up to at most one.

Also the Bernoulli case can be studied in the view of Algorithmic
Information Theory. We call this the \emph{universal setup}: Given
a universal Turing machine, the related class of Bernoulli
distributions is isomorphic to the countable set of computable
reals in $[0,1]$. The description length $\K w(\th)$ of a
parameter $\th\in[0,1]$ is then given by the length of its
shortest program. A prior weight may then be defined by $2^{-\K
w(\th)}$. (If a string $x=x_1x_2\ldots x_{t-1}$ is generated by a
Bernoulli distribution with computable parameter
$\th_0\in[0,1]$, then with high probability the two-part
complexity of $x$ with respect to the Bernoulli class does not
exceed its algorithmic complexity by more than a constant, as
shown by Vovk \cite{Vovk:97}. That is, the two-part complexity
with respect to the Bernoulli class {\it is} the shortest
description, save for an additive constant.)

Many Machine Learning tasks are or can be reduced to sequence
prediction tasks. An important example is classification. The task
of classifying a new instance $z_n$ after having seen
(instance,class) pairs $(z_1,c_1),...,(z_{n-1},c_{n-1})$ can be
phrased as to predict the continuation of the sequence
$z_1c_1...z_{n-1}c_{n-1}z_n$. Typically the (instance,class) pairs
are i.i.d. Cumulative loss bounds for prediction usually
generalize to prediction \emph{conditionalized} to some inputs
\cite{Poland:05mdlreg}. Then we can solve classification problems
in the standard form. It is not obvious if and how the proofs in
this paper can be conditionalized.

Our main tool for obtaining results is the Kullback-Leibler
divergence. Lemmata for this quantity are stated in Section
\ref{secKL}. Section \ref{secLB} shows that the exponential error
bound obtained in \cite{Poland:04mdl} is sharp in general. In
Section \ref{secUB}, we give an upper bound on the instantaneous
and the cumulative losses. The latter bound is small e.g.\ under
certain conditions on the distribution of the weights, this is the
subject of Section \ref{secUDW}. Section \ref{secUC} treats the
universal setup. Finally, in Section \ref{secDC} we discuss the
results and give conclusions.

\section{Kullback-Leibler Divergence}\label{secKL}

Let $\BBB=\{0,1\}$ and consider finite strings $x\in\BBB^*$ as
well as infinite sequences $x\ltinf\in\BBB^\infty$, with the first
$n$ bits denoted by $x_{1:n}$. If we know
that $x$ is generated by an i.i.d random variable, then
$P(x_i=1)=\th_0$ for all $1\leq i\leq\l(x)$ where $\l(x)$ is the
length of $x$. Then $x$ is called a Bernoulli sequence, and
$\th_0\in\Theta\subset[0,1]$ the {\it true parameter}. In the
following we will consider only countable $\Theta$, e.g.\ the set
of all computable numbers in $[0,1]$.

Associated with each $\th\in\Theta$, there is a {\it complexity}
or description length $\K w(\th)$ and a {\it weight} or
(semi)probability $w_\th=2^{-\K w(\th)}$. The complexity will
often but need not be a natural number. Typically, one assumes
that the weights sum up to at most one,
$\sum_{\th\in\Theta}w_\th\leq 1$. Then, by the Kraft inequality,
for all $\th\in\Theta$ there exists a prefix-code of length
$\K w(\th)$. Because of this correspondence, it is only a matter of
convenience whether results are developed in terms of description
lengths or probabilities. We will choose the former way. We won't
even need the condition $\sum_\th w_\th\leq 1$ for most of the
following results. This only means that $\K w$ cannot be
interpreted as a prefix code length, but does not cause other
problems.

Given a set of distributions $\Theta\subset [0,1]$, complexities
$\big(\K w(\th)\big)_{\th\in\Theta}$, a true distribution
$\th_0\in\Theta$, and some observed string $x\in\BBB^*$, we
define an {\it MDL estimator}%
\footnote{\label{footnoteMDL}Precisely, we define a MAP (maximum a posteriori)
estimator. For two reasons, our definition might not be considered
as MDL in the strict sense. First, MDL is often associated with a
specific prior, while we admit arbitrary priors. Second and more
importantly, when coding some data
$x$, one can exploit the fact that once the parameter $\th^x$ is
specified, only data which leads to this $\th^x$ needs to be
considered. This allows for a description shorter than $\K
w(\th^x)$. Nevertheless, the
\emph{construction principle} is commonly termed MDL, compare
e.g.\ the ``ideal MDL" in \cite{Vitanyi:00}.}:
\beqn
  \th^x =\arg\max_{\th\in\Theta}\{ w_\th P(x|\th)\}.
\eeqn
Here, $P(x|\th)$ is the probability of observing $x$ if
$\th$ is the true parameter. Clearly,
$P(x|\th)=\th^{\eins(x)}(1-\th)^{\l(x)-\eins(x)}$,
where $\eins(x)$ is the number of ones in $x$. Hence
$P(x|\th)$ depends only on $\l(x)$ and
$\eins(x)$. We therefore see
\bqa\label{eq:argmax}\label{defthhat}
  \th^x\ =\ \th^{(\alpha,n)}
  & = & \arg\max_{\th\in\Theta}\{ w_\th
  \left(\th^{\alpha}(1-\th)^{1-\alpha}\right)^n\}\\
  \nonumber
  & = & \arg\min_{\th\in\Theta}\{n\!\cdot\!D(\alpha\|\th)+\K w(\th)\!\cdot\ln 2\},
\eqa
where $n=\l(x)$ and $\alpha:=\frac{\eins(x)}{\l(x)}$ is the {\it
observed fraction} of ones and
$$\textstyle D(\alpha\|\th)=\alpha\ln\frac\alpha\th+(1-\alpha)\ln\frac{1-\alpha}{1-\th}$$
is the Kullback-Leibler divergence. The second line of
(\ref{defthhat}) also explains the name MDL, since we choose the
$\th$ which minimizes the joint description of model $\th$ and the
data $x$ given the model.

We also define the {\it
extended Kullback-Leibler divergence}
\beq
\label{eq:ExtendedKullbackLeibler}
  D^\alpha(\th\|\tilde\th)
  = \alpha\ln\frac{\th}{\tilde\th}+(1-\alpha)\ln\frac{1-\th}{1-\tilde\th}
  = D(\alpha\|\tilde\th) - D(\alpha\|\th).
\eeq
It is easy to see that $D^\alpha(\th\|\tilde\th)$ is linear in $\alpha$,
$D^\th(\th\|\tilde\th)=D(\th\|\tilde\th)$ and
$D^{\tilde\th}(\th\|\tilde\th)=-D(\tilde\th\|\th)$,
and $\frac{d}{d\alpha}D^\alpha(\th\|\tilde\th)>0$
iff $\th>\tilde\th$. Note that $D^\alpha(\th\|\tilde\th)$
may be also defined for the general i.i.d.\ case, i.e.\ if
the alphabet has more than two symbols.

Let $\th,\tilde\th\in\Theta$ be two parameters, then it follows
from (\ref{eq:argmax}) that in the process of choosing the MDL
estimator, $\th$ is being preferred to $\tilde\th$ iff
\beq
\label{eq:Beating}
nD^\alpha(\th\|\tilde\th)\geq \ln 2 \cdot \big(\K w(\th)-\K
w(\tilde\th)\big)
\eeq
with $n$ and $\alpha$ as before. We also say that then $\th$ {\it
beats} $\tilde\th$. It is immediate that for increasing $n$ the
influence of the complexities on the selection of the maximizing
element decreases.
We are now interested in the {\it total expected square prediction
error} (or cumulative square loss) of the MDL estimator
\beqn
\sum_{n=1}^\infty \Expect(\th^{x\leqn}-\th_0)^2.
\eeqn
In terms of \cite{Poland:04mdl}, this is the {\it static MDL
prediction} loss, which means that a predictor/estimator $\th^x$
is chosen according to the current observation $x$. (As already
mentioned, the terms predictor and estimator coincide for static
MDL and Bernoulli classes.) The \emph{dynamic} method on the other
hand would consider both possible continuations $x0$ and $x1$ and
predict according to $\th^{x0}$ \emph{and} $\th^{x1}$. In the
following, we concentrate on static predictions. They are also
preferred in practice, since computing only one model is more
efficient.

Let $A_n=\big\{\frac{k}{n}:0\leq k\leq n\big\}$. Given the true
parameter $\th_0$ and some $n\in\NNN$, the {\it expectation} of a
function $f^{(n)}:\{0,\ldots,n\}\to\RRR$ is given by
\beq\label{defE}
  \Expect f^{(n)}= \sum_{\alpha\in A_n} p(\alpha|n) f(\alpha n),
  \mbox{ where }
  p(\alpha|n)={n \choose k}\Big(\th_0^\alpha (1-\th_0)^{1-\alpha}\Big)^n.
\eeq
(Note that the probability $p(\alpha|n)$ depends on $\th_0$, which
we do not make explicit in our notation.) Therefore,
\beq
  \sum_{n=1}^\infty \Expect(\th^{x\leqn}-\th_0)^2=
  \sum_{n=1}^\infty \sum_{\alpha\in A_n} p(\alpha|n) (\th^{(\alpha,n)}-\th_0)^2,
\eeq
Denote the relation
$f=O(g)$ by $f\leqm g$. Analogously define ``$\geqm$"
and ``$\eqm$". From \cite[Corollary 12]{Poland:04mdl}, we
immediately obtain the following result.

\begin{Theorem}
\label{th:Previous}
The cumulative loss bound
$\sum_n \Expect(\th^{x\leqn}-\th_0)^2\leqm 2^{\K w(\th_0)}$ holds.
\end{Theorem}

This is the ``slow" convergence result mentioned in the
introduction. In contrast, for a Bayes mixture, the total expected
error is bounded by $\K w(\th_0)$ rather than $2^{\K w(\th_0)}$
(see \cite{Solomonoff:78} or \cite[Th.1]{Hutter:01alpha}). An
upper bound on $\sum_n \Expect(\th^{x\leqn}-\th_0)^2$ is termed as
\emph{convergence in mean sum} and implies convergence
$\th^{x\leqn}\to\th_0$ with probability 1 (since otherwise the sum
would be infinite).

We now establish relations between the Kullback-Leibler divergence
and the quadratic distance. We call bounds of this type {\it
entropy inequalities}.

\begin{Lemma} \label{Lemma:EntropyIneq}
Let $\th,\tilde\th\in(0,1)$. Let
$\th^*=\arg\min\{|\th-\frac{1}{2}|,|\tilde\th-\frac{1}{2}|\}$, i.e.\ $\th^*$
is the element from $\{\th,\tilde\th\}$ which is closer to $\frac{1}{2}$.
Then the following assertions hold.
\bqan
(i) & D(\th\|\tilde\th) \ \geq \ 2\cdot(\th-\tilde\th)^2
&\ \forall\ \th,\tilde\th\in(0,1),\\
(ii) & D(\th\|\tilde\th) \ \leq \ \mbox{$\frac{8}{3}$}(\th-\tilde\th)^2
& \wenn \th,\tilde\th\in[\mbox{$\frac{1}{4}$},\mbox{$\frac{3}{4}$}],\\
(iii) & D(\th\|\tilde\th) \ \geq \ \frac{(\th-\tilde\th)^2}{2\th^*(1-\th^*)}
& \wenn \th,\tilde\th\leq\mbox{$\frac{1}{2}$},\\
(iv) & D(\th\|\tilde\th) \ \leq \ \frac{3(\th-\tilde\th)^2}{2\th^*(1-\th^*)}
& \wenn \th\leq\mbox{$\frac{1}{4}$} \und \tilde\th\in[\mbox{$\frac{\th}{3}$},3\th],\\
(v) & D(\tilde\th\|\th) \ \geq \ \tilde\th(\ln\tilde\th-\ln\th-1)
&\ \forall\ \th,\tilde\th\in(0,1),\\
(vi) & D(\th\|\tilde\th) \ \leq \ \mbox{$\frac{1}{2}$}\tilde\th
& \wenn \th\leq\tilde\th\leq\mbox{$\frac{1}{2}$},\\
(vii) & D(\th\|\th\cdot 2^{-j}) \ \leq \ j\cdot\th
& \wenn \th\leq\mbox{$\frac{1}{2}$}\und j\geq  1,\\
(viii) & D(\th\|1-2^{-j}) \ \leq \ j
& \wenn \th\leq\mbox{$\frac{1}{2}$} \und j\geq 1.
\eqan
Statements $(iii)-(viii)$ have symmetric counterparts
for $\th\geq\frac{1}{2}$.
\end{Lemma}

The first two statements give upper and lower bounds for the
Kullback-Leibler divergence in terms of the quadratic distance.
They express the fact that the Kullback-Leibler divergence is
locally quadratic. So do the next two statements, they will be
applied in particular if $\th$ is located close to the boundary of
$[0,1]$. Statements $(v)$ and $(vi)$ give bounds in terms of the
absolute distance, i.e.\ ``linear" instead of quadratic. They are
mainly used if $\tilde\th$ is relatively far from $\th$. Note that
in $(v)$, the position of $\th$ and $\tilde\th$ are inverted. The
last two inequalities finally describe the behavior of the
Kullback-Leibler divergence as its second argument tends to the
boundary of $[0,1]$. Observe that this is logarithmic in the
inverse distance to the boundary.

\begin{Proof}
$(i)$ This is standard, see e.g. \cite{Li:97}. It is shown
similarly as $(iii)$.

$(ii)$ Let $f(\eta)=D(\th\|\eta)-\frac{8}{3}(\eta-\th)^2$, then we show
$f(\eta)\leq 0$ for $\eta\in[\frac{1}{4},\frac{3}{4}]$. We have
that $f(\th)=0$ and
\beqn
f'(\eta)=\frac{\eta-\th}{\eta(1-\eta)}-\frac{16}{3}(\eta-\th).
\eeqn
This difference is nonnegative if and only $\eta-\th\leq 0$ since
$\eta(1-\eta)\geq\frac{3}{16}$. This implies $f(\eta)\leq 0$.

$(iii)$
Consider the function
\beqn
f(\eta)=D(\th\|\eta)-\frac{(\th-\eta)^2}{2\max\{\th,\eta\}(1-\max\{\th,\eta\})}.
\eeqn
We have to show that $f(\eta)\geq 0$ for all $\eta\in(0,\frac{1}{2}]$. It is obvious
that $f(\th)=0$. For $\eta\leq\th$,
\beqn
f'(\eta)=\frac{\eta-\th}{\eta(1-\eta)}-\frac{\eta-\th}{\th(1-\th)}\leq 0
\eeqn
holds since $\eta-\th\leq 0$ and $\th(1-\th)\geq \eta(1-\eta)$.
Thus, $f(\eta)\geq 0$ must be valid for $\eta\leq\th$.
On the other hand if $\eta\geq\th$, then
\beqn
f'(\eta)=\frac{\eta-\th}{\eta(1-\eta)}-\left[\frac{\eta-\th}{\eta(1-\eta)}
-\frac{(\eta-\th)^2(1-2\eta)}{2\eta^2(1-\eta)^2}\right]\geq 0
\eeqn
is true. Thus $f(\eta)\geq 0$ holds in this case, too.

$(iv)$
We show that
\beqn
f(\eta)=D(\th\|\eta)-\frac{3(\th-\eta)^2}{2\max\{\th,\eta\}(1-\max\{\th,\eta\})}\leq 0
\eeqn
for $\eta\in[\frac{\th}{3},3\th]$.
If $\eta\leq\th$, then
\beqn
f'(\eta)=\frac{\eta-\th}{\eta(1-\eta)}-\frac{3(\eta-\th)}{\th(1-\th)}\geq 0
\eeqn
since $3\eta(1-\eta)\geq\th(1-\eta)\geq\th(1-\th)$.
If $\eta\geq\th$, then
\beqn
f'(\eta)=\frac{\eta-\th}{\eta(1-\eta)}-3\cdot \left[\frac{\eta-\th}{\eta(1-\eta)}
-\frac{(\eta-\th)^2(1-2\eta)}{2\eta^2(1-\eta)^2}\right]\leq 0
\eeqn
is equivalent to $4\eta(1-\eta)\geq 3(\eta-\th)(1-2\eta)$, which is fulfilled
if $\th\leq\frac{1}{4}$ and $\eta\leq 3\th$ as an elementary computation
verifies.

$(v)$ Using $-\ln(1\!-\!u)\leq\frac{u}{1-u}$, one obtains
\bqan
D(\tilde\th\|\th) & = &
  \tilde\th\ln\frac{\tilde\th}{\th}+(1-\tilde\th)\ln\frac{1-\tilde\th}{1-\th}
\ \geq \
  \tilde\th\ln\frac{\tilde\th}{\th}+(1-\tilde\th)\ln(1-\tilde\th)
\\ & \geq &
  \tilde\th\ln\frac{\tilde\th}{\th}-\tilde\th
\ = \
  \tilde\th(\ln\tilde\th-\ln\th-1)
\eqan

$(vi)$
This follows from
$D(\th\|\tilde\th)\leq -\ln(1-\tilde\th)\leq\frac{\tilde\th}{1-\tilde\th}\leq
\frac{\tilde\th}{2}$. The last two statements
$(vii)$ and $(viii)$ are even easier.
\end{Proof}

In the above entropy inequalities we have left out the extreme cases
$\th,\tilde\th\in\{0,1\}$. This is for simplicity and convenience
only. Inequalities $(i)-(iv)$ remain valid for
$\th,\tilde\th\in\{0,1\}$ if the fraction $\frac{0}{0}$ is
properly defined. However, since the extreme cases will need to be
considered separately anyway, there is no requirement for the
extension of the lemma. We won't need $(vi)$ and $(viii)$ of Lemma
\ref{Lemma:EntropyIneq} in the sequel.

We want to point out that although we have proven Lemma
\ref{Lemma:EntropyIneq} only for the case of binary alphabet,
generalizations to arbitrary alphabet are likely to hold. In fact,
$(i)$ does hold for arbitrary alphabet, as shown in
\cite{Hutter:01alpha}.

It is a well-known fact that the binomial distribution may be
approximated by a Gaussian. Our next goal is to establish upper
and lower bounds for the binomial distribution. Again we leave out
the extreme cases.

\begin{Lemma} \label{Lemma:BinomialBounds}
Let $\th_0\in(0,1)$ be the true parameter, $n\geq 2$ and $1\leq
k\leq n-1$, and $\alpha=\frac{k}{n}$. Then the following
assertions hold.
\bqan
(i) && p(\alpha|n)\leq
\frac{1}{\sqrt{2\pi\alpha(1-\alpha)n}}\exp\big(-nD(\alpha\|\th_0)\big),\\
(ii) && p(\alpha|n)\geq
\frac{1}{\sqrt{8\alpha(1-\alpha)n}}\exp\big(-nD(\alpha\|\th_0)\big).
\eqan
\end{Lemma}
The lemma gives a
quantitative assertion about the Gaussian approximation to a
binomial distribution. The upper bound is sharp for $n\to\infty$
and fixed $\alpha$. Lemma \ref{Lemma:BinomialBounds} can be easily
combined with Lemma \ref{Lemma:EntropyIneq}, yielding Gaussian
estimates for the Binomial distribution.

\begin{Proof}
Stirling's formula is a well-known result from calculus. In a refined version,
it states that for any $n\geq 1$ the factorial $n!$ can be bounded from below
and above by
\beqn
\sqrt{2\pi n}\cdot n^n\exp\left(-n+\frac{1}{12n+1}\right)\leq n!\leq
\sqrt{2\pi n}\cdot n^n\exp\left(-n+\frac{1}{12n}\right).
\eeqn
Hence,
\bqan
  p(\alpha,n) & = & \frac{n!}{k!(n-k)!}\th_0^k(1-\th_0)^{n-k}
\\
  & \leq & \frac{\sqrt{n}\cdot n^n\exp\left(\frac{1}{12n}\right)\th_0^k(1-\th_0)^{n-k}}
  {\sqrt{2\pi k(n-k)}\cdot k^k(n-k)^{n-k}\exp\left(\frac{1}{12k+1}+\frac{1}{12(n-k)+1}\right)}
\\
  & = & \frac{1}{\sqrt{2\pi\alpha(1-\alpha)n}}\exp\left(
  -n\cdot D(\alpha\|\th_0)+\frac{1}{12n}-\frac{1}{12k+1}-\frac{1}{12(n-k)+1} \right)
\\
  & \leq & \frac{1}{\sqrt{2\pi\alpha(1-\alpha)n}}\exp\big(-nD(\alpha\|\th_0)\big).
\eqan
The last inequality is valid since
$\frac{1}{12n}-\frac{1}{12k+1}-\frac{1}{12(n-k)+1}<0$ for all $n$ and $k$,
which is easily verified using elementary computations. This
establishes $(i)$.

In order to show $(ii)$, we observe
\bqan p(\alpha,n) & \geq &
  \frac{1}{\sqrt{2\pi\alpha(1-\alpha)n}}\exp\left(
  -n\cdot D(\alpha\|\th_0)+\frac{1}{12n+1}-\frac{1}{12k}-\frac{1}{12(n-k)} \right)
\\
  & \geq & \frac{\exp({1\over 37}-{1\over 8})}{\sqrt{2\pi\alpha(1-\alpha)n}}\exp\big(-nD(\alpha\|\th_0)\big)
  \mbox{\quad for\quad} n\geq 3.
\eqan
Here the last inequality follows from the fact that
$\frac{1}{12n+1}-\frac{1}{12k}-\frac{1}{12(n-k)}$
is minimized for $n=3$ (and $k=1$ or $2$), if we exclude $n=2$,
and $\exp({1\over 37}-{1\over 8})\geq\sqrt\pi/2$. For $n=2$ a
direct computation establishes the lower bound.
\end{Proof}

\begin{Lemma} \label{Lemma:IntegralEstimate}
Let $z\in\RRR^+$, then
\bqan
  (i) && \frac{\sqrt\pi}{2z^3}-\frac{1}{z\sqrt{2e}}
  \leq \sum_{n=1}^\infty \sqrt{n}\cdot\exp(-z^2n)
  \leq \frac{\sqrt\pi}{2z^3}+\frac{1}{z\sqrt{2e}} \und\\
  (ii) && \sum_{n=1}^\infty n^{-\frac{1}{2}}\exp(-z^2n) \leq
  \sqrt\pi/z.
\eqan
\end{Lemma}

\begin{Proof}
$(i)$ The function $f(u)=\sqrt{u}\exp(-z^2u)$ increases for $u\leq\frac{1}{2z^2}$
and decreases for $u\geq\frac{1}{2z^2}$. Let
$N=\max\{n\in\NNN:f(n)\geq f(n-1)\}$, then it is easy to see that
\bqan
\sum_{n=1}^{N-1} f(n) & \leq & \int_0^N f(u)\ du
\ \leq \ \sum_{n=1}^N f(n) \und\\
\sum_{n=N+1}^\infty f(n) & \leq & \int_N^\infty f(u)\ du
\ \leq \ \sum_{n=N}^\infty f(n) \mbox{ and thus}\\
\sum_{n=1}^\infty f(n)-f(N) & \leq & \int_0^\infty f(u)\ du
\ \leq \ \sum_{n=1}^\infty f(n)+f(N)
\eqan
holds. Moreover, $f$ is the derivative of the function
\beqn
F(u)=-\frac{\sqrt u\exp(-z^2 u)}{z^2}+\frac{1}{z^3}\int_0^{z\sqrt u}\exp(-v^2)\ dv.
\eeqn
Observe $f(N)\leq
f(\frac{1}{2z^2})=\frac{\exp(-\frac{1}{2})}{z\cdot\sqrt 2}$ and
$\int_0^\infty\exp(-v^2)dv=\frac{\sqrt\pi}{2}$ to obtain the
assertion.

$(ii)$ The function $f(u)=u^{-\frac{1}{2}}\exp(-z^2u)$
decreases monotonically on $(0,\infty)$ and is the derivative of
$F(u)=2z^{-1}\int_0^{z\sqrt{u}}\exp(-v^2)dv$. Therefore,
\beqn
\sum_{n=1}^\infty f(n) \leq
\int_0^\infty f(u)\ du = \sqrt\pi/z
\eeqn
holds.
\end{Proof}

\section{Lower Bound}\label{secLB}

We are now in the position to prove that even for Bernoulli
classes the upper bound from Theorem \ref{th:Previous} is sharp in
general.

\begin{Prop} \label{Prop:Counterex}
Let $\th_0=\frac{1}{2}$ be the true parameter generating sequences
of fair coin flips. Assume $\Theta=\{\th_0,\th_1,\ldots,\th_{2^N-1}\}$ where $\th_k=\frac{1}{2}+2^{-k-1}$ for $k\geq 1$.
Let all complexities be equal, i.e.\ $\K w(\th_0)=\ldots=\K w(\th_{2^N-1})=N$. Then
\beqn
  \sum_{n=1}^\infty \Expect (\th_0-\th^x)^2\geq \mbox{$\frac{1}{84}$}\big(2^N-5\big)
  \eqm 2^{\K w(\th_0)}.
\eeqn
\end{Prop}

\begin{Proof}
Recall that $\th^x=\th^{(\alpha,n)}$ the maximizing element for
some observed sequence $x$ only depends on the length $n$ and the
observed fraction of ones $\alpha$. In order to obtain an estimate
for the total prediction error $\sum_n\Expect (\th_0-\th^x)^2$,
partition the interval $[0,1]$ into $2^N$ disjoint intervals
$I_k$, such that $\bigcup_{k=0}^{2^N-1} I_k=[0,1]$.
Then consider the contributions for the observed fraction $\alpha$
falling in $I_k$ {\it separately}:
\beq
\label{eq:Ck}
C(k)= \sum_{n=1}^\infty \sum_{\alpha\in A_n\cap
I_k}p(\alpha|n)(\th^{(\alpha,n)}-\th_0)^2
\eeq
(compare (\ref{defE})). Clearly, $\sum_n\Expect
(\th_0-\th^x)^2=\sum_k C(k)$ holds. We define the partitioning
$(I_k)$ as
$I_0=[0,\frac{1}{2}+2^{-2^N})=[0,\th_{2^N-1})$,
$I_1=[\frac{3}{4},1]=[\th_1,1]$, and
\beqn I_k=[\th_k,\th_{k-1}) \for_all 2\leq k\leq 2^N-1.
\eeqn
Fix $k\in\{2,\ldots,2^N-1\}$ and assume $\alpha\in I_k$. Then
\beqn
\th^{(\alpha,n)}=\arg\min_\th\{nD(\alpha\|\th)+\K w(\th)\ln 2\}
=\arg\min_\th\{nD(\alpha\|\th)\}\in\{\th_k,\th_{k-1}\}
\eeqn
according to (\ref{eq:argmax}). So clearly
$(\th^{(\alpha,n)}-\th_0)^2\geq(\th_k-\th_0)^2=2^{-2k-2}$ holds. Since
$p(\alpha|n)$ decreases for increasing $|\alpha-\th_0|$, we have
$p(\alpha|n)\geq p(\th_{k-1}|n)$. The interval $I_k$ has length $2^{-k-1}$,
so there are at least $\lfloor n 2^{-k-1}\rfloor\geq n 2^{-k-1}-1$
observed fractions $\alpha$ falling in the interval. From
(\ref{eq:Ck}), the total contribution of $\alpha\in I_k$ can be
estimated by
\beqn
C(k)\geq\sum_{n=1}^\infty 2^{-2k-2}(n2^{-k-1}-1)p(\th_{k-1}|n).
\eeqn
Note that the terms in the sum even become negative for small $n$,
which does not cause any problems. We proceed with
\beqn
p(\th_{k-1}|n) \geq
\frac{1}{\sqrt{8\cdot 2^{-2}n}}\exp
\big[-nD\big(\mbox{$\frac{1}{2}$}+2^{-k}\|\mbox{$\frac{1}{2}$}\big)\big]
\geq
\frac{1}{\sqrt{2n}}\exp
\big[-n\mbox{$\frac{8}{3}$} 2^{-2k}\big]
\eeqn
according to Lemma \ref{Lemma:BinomialBounds} and Lemma
\ref{Lemma:EntropyIneq} (ii). By Lemma
\ref{Lemma:IntegralEstimate} $(i)$ and $(ii)$, we have
\bqan
\sum_{n=1}^\infty \sqrt n\exp \big[-n\mbox{$\frac{8}{3}$} 2^{-2k}\big]
& \geq & \frac{\sqrt\pi}{2}\left(\frac{3}{8}\right)^{\frac{3}{2}}2^{3k}-
\frac{1}{\sqrt{2e}}\sqrt{\frac{3}{8}} 2^{k}\und\\
-\sum_{n=1}^\infty n^{-\frac{1}{2}}\exp \big[-n\mbox{$\frac{8}{3}$} 2^{-2k}\big]
& \geq & -\sqrt\pi \sqrt{\frac{3}{8}} 2^{k}.
\eqan
Considering only $k\geq 5$, we thus obtain
\bqan
C(k)& \geq &
\frac{1}{\sqrt{2}}\sqrt{\frac{3}{8}}2^{-2k-2}
\left[\frac{3\sqrt\pi}{16}2^{2k-1}-\frac{1}{\sqrt{2e}} 2^{-1}-\sqrt\pi 2^{k}
\right]
\\ & \geq &
\frac{\sqrt 3}{16}
\left[3\sqrt\pi 2^{-5}-\frac{1}{\sqrt{2e}} 2^{-2k-1}-\sqrt\pi 2^{-k}
\right]
 \geq
\frac{\sqrt{3\pi}}{8}2^{-5}-\frac{\sqrt 3}{16\sqrt{2e}}2^{-11}
>\frac{1}{84}.
\eqan
Ignoring the contributions for $k\leq 4$, this implies the
assertion.
\end{Proof}

This result shows that if the parameters and their weights are
chosen in an appropriate way, then the total expected error is of
order $w_0^{-1}$ instead of $\ln w_0^{-1}$. Interestingly, this
outcome seems to depend on the arrangement and the weights of the
{\it false} parameters rather than on the weight of the {\it true}
one. One can check with moderate effort that the proposition still
remains valid if e.g.\ $w_0$ is twice as large as the other
weights. Actually, the proof of Proposition \ref{Prop:Counterex}
shows even a slightly more general result, namely admitting
additional arbitrary parameters with larger complexities:

\begin{Cor} \label{Cor:Counterex}
Let $\Theta=\{\th_k:k\geq 0\}$, $\th_0=\frac{1}{2}$,
$\th_k=\frac{1}{2}+2^{-k-1}$ for $1\leq k\leq 2^N-2$,
and $\th_k\in[0,1]$ arbitrary for $k\geq 2^N-1$. Let
$\K w(\th_k)=N$ for $0\leq k\leq 2^N-2$ and $\K w(\th_k)>N$ for $k\geq 2^N-1$.
Then $\sum_n \Expect (\th_0-\th^x)^2\geq \frac{1}{84}(2^N-6)$
holds.
\end{Cor}

We will use this result only for Example \ref{Ex:Sensitive}. Other
and more general assertions can be proven similarly.

\section{Upper Bounds}\label{secUB}

Although the cumulative error may be large, as seen in the
previous section, the instantaneous error is always small. It is
easy to demonstrate this for the Bernoulli case, to which we
restrict in this paper. Much more general results have been
obtained for arbitrary classes of i.i.d.\ models
\cite{Li:99,Zhang:04}. Strong instantaneous bounds hold in
particular if MDL is modified by replacing the factor $\ln 2$ in
(\ref{defthhat}) by something larger (e.g. $(1+\eps)\ln 2$) such
that complexity is penalized slightly more than usually. Note that
our cumulative bounds are incomparable to these and other
instantaneous bounds.

\begin{Prop}\label{propSUB}
For $n\geq 3$, the expected instantaneous square loss is bounded as follows:
\beqn
  \Expect (\th_0-\hat\th^{x_{1:n}})^2
    \leq \frac{(\ln 2)\K w(\th_0)}{2n}+\frac{\sqrt{2(\ln 2)\K w(\th_0)\ln n}}{n}+
    \frac{6\ln n}{n}.
\eeqn
\end{Prop}

\begin{Proof}
We give an elementary proof for the case
$\th_0\in(\frac{1}{4},\frac{3}{4})$ only. Like in the proof of
Proposition \ref{Prop:Counterex}, we consider the contributions of
different $\alpha$ separately. By Hoeffding's inequality,
$\Prob(|\alpha-\th_0|\geq\frac{c}{\sqrt n})\leq2e^{-2c^2}$
for any $c>0$. Letting $c=\sqrt{\ln n}$, the contributions by
these $\alpha$ are thus bounded by $\frac{2}{n^2}\leq\frac{\ln
n}{n}$.

On the other hand, for $|\alpha-\th_0|\leq\frac{c}{\sqrt n}$,
recall that $\th_0$ \emph{beats} any $\th$ iff (\ref{eq:Beating})
holds. According to $\K w(\th)\geq 0$,
$|\alpha-\th_0|\leq\frac{c}{\sqrt n}$, and Lemma
\ref{Lemma:EntropyIneq} $(i)$ and $(ii)$, (\ref{eq:Beating}) is
already implied by
$
|\alpha-\th|\geq\sqrt{\frac{\frac{1}{2}(\ln 2)\K
w(\th_0)+\frac{4}{3}c^2}{n}}.
$
Clearly, a contribution only occurs if $\th$ beats
$\th_0$, therefore if the opposite inequality holds.
Using $|\alpha-\th_0|\leq\frac{c}{\sqrt n}$ again and the triangle
inequality, we obtain that
\beqn(\th-\th_0)^2\leq\frac{5c^2+
\frac{1}{2}(\ln 2)\K w(\th_0)+\sqrt{2(\ln 2)\K w(\th_0)c^2}}{n}
\eeqn
in this case. Since we have chosen $c=\sqrt{\ln n}$, this implies
the assertion.
\end{Proof}

One can improve the bound in Proposition \ref{propSUB} to
$\Expect (\th_0-\hat\th^{x_{1:n}})^2\leqm \frac{\K w(\th_0)}{n}$
by a refined argument, compare \cite{Barron:91}. But the
high-level assertion is the same: Even if the cumulative upper
bound may be infinite, the instantaneous error converges
rapidly to 0. Moreover, the convergence speed depends on
$\K w(\th_0)$ as opposed to $2^{\K w(\th_0)}$. Thus $\hat\th$
tends to $\th_0$ rapidly in probability (recall that the assertion
is not strong enough to conclude almost sure convergence). The
proof does not exploit $\sum w_\th\leq 1$, but only
$w_\th\leq 1$, hence the assertion even holds for a maximum
likelihood estimator (i.e.\ $w_\th=1$ for all
$\th\in\Theta$). The theorem generalizes to i.i.d.\ classes. For
the example in Proposition \ref{Prop:Counterex}, the instantaneous
bound implies that the bulk of losses occurs very late. This does
\emph{not} hold for general (non-i.i.d.) model classes:
The total loss up to time $n$ in \cite[Example 9]{Poland:04mdl}
grows linearly in $n$.

We will now state our main positive result that upper bounds the
cumulative loss in terms of the negative logarithm of the true
weight and the {\it arrangement} of the false parameters. The
proof is similar to that of Proposition \ref{Prop:Counterex}. We
will only give the proof idea here and defer the lengthy and tedious
technical details to the appendix.

Consider the cumulated sum square error $\sum_n \Expect
\big(\th^{(\alpha,n)}-\th_0\big)^2$. In order to upper bound this
quantity, we will partition the open unit interval $(0,1)$ into a
sequence of intervals $(I_k)_{k=1}^\infty$, each of measure
$2^{-k}$. (More precisely: Each $I_k$ is either an interval or a
union of two intervals.) Then we will estimate the contribution of
each interval to the cumulated square error,
\beqn
C(k)=\sum_{n=1}^\infty \sum_{\alpha\in A_n,\th^{(\alpha,n)}\in
I_k} p(\alpha|n) (\th^{(\alpha,n)}-\th_0)^2
\eeqn
(compare (\ref{defE}) and (\ref{eq:Ck})). Note that
$\th^{(\alpha,n)}\in I_k$ precisely reads $\th^{(\alpha,n)}\in
I_k\cap\Theta$, but for convenience we generally assume
$\th\in\Theta$ for all $\th$ being considered. This
partitioning is also used for $\alpha$, i.e. define the
contribution $C(k,j)$ of $\th\in I_k$ where $\alpha\in I_j$ as
\beqn
C(k,j)=\sum_{n=1}^\infty \sum_{\alpha\in A_n\cap
I_j,\th^{(\alpha,n)}\in I_k} p(\alpha|n)
(\th^{(\alpha,n)}-\th_0)^2.
\eeqn
We need to distinguish between
$\alpha$ that are located close to $\th_0$ and
$\alpha$ that are located far from $\th_0$. ``Close" will
be roughly equivalent to $j>k$, ``far" will be
approximately $j\leq k$.
So we get
$\sum_n \Expect \big(\th^{(\alpha,n)}-\th_0\big)^2
= \sum_{k=1}^\infty C(k) = \sum_k\sum_j C(k,j)$. In the proof,
\beqn
p(\alpha|n)\leqm
\big[n\alpha(1-\alpha)\big]^{-\frac{1}{2}}\exp\big[-nD(\alpha\|\th_0)\big]
\eeqn
is often applied, which holds by Lemma \ref{Lemma:BinomialBounds}
(recall that $f\leqm g$ stands for
$f=O(g)$). Terms like $D(\alpha\|\th_0)$, arising in
this context and others, can be further estimated using Lemma
\ref{Lemma:EntropyIneq}. We now give the constructions of
intervals $I_k$ and complementary intervals $J_k$.

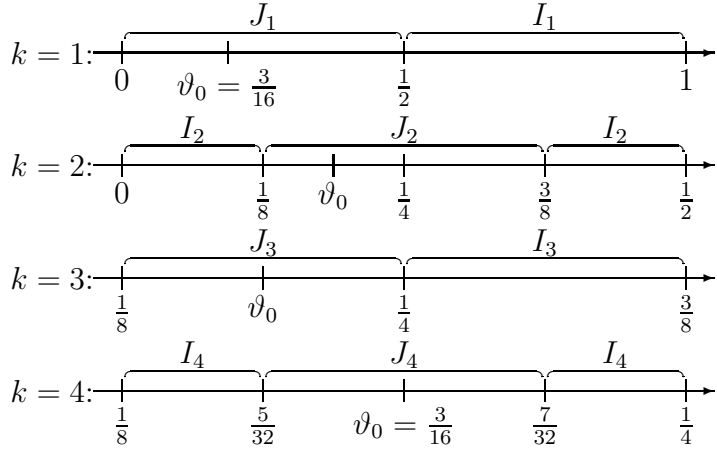
\begin{figure}[t]\begin{center}
\setlength{\unitlength}{7.5cm}
\begin{picture}(1.3,0.8)
\multiput(0.15,0.2)(0,0.2){4}{\vector(1,0){1.1}}
\multiput(0.2,0.18)(0,0.2){4}{\line(0,1){0.04}}
\multiput(1.2,0.18)(0,0.2){4}{\line(0,1){0.04}}
\put(0.2,0.77){\makebox(0,0)[t]{$0$}}
\put(0.2,0.57){\makebox(0,0)[t]{$0$}}
\put(1.2,0.77){\makebox(0,0)[t]{$1$}}
\put(1.2,0.57){\makebox(0,0)[t]{{$\frac{1}{2}$}}}
\multiput(0.7,0.18)(0,0.2){4}{\line(0,1){0.04}}
\put(0.3875,0.78){\line(0,1){0.04}}
\put(0.3875,0.77){\makebox(0,0)[t]{$\th_0=\frac{3}{16}$}}
\put(0.7,0.77){\makebox(0,0)[t]{$\frac{1}{2}$}}
\put(0.575,0.58){\line(0,1){0.04}}
\put(0.575,0.57){\makebox(0,0)[t]{$\th_0$}}
\put(0.7,0.57){\makebox(0,0)[t]{$\frac{1}{4}$}}
\put(0.45,0.825){\oval(0.49,0.02)[t]}
\put(0.95,0.825){\oval(0.49,0.02)[t]}
\put(0.45,0.84){\makebox(0,0)[b]{$J_1$}}
\put(0.95,0.84){\makebox(0,0)[b]{$I_1$}}

\put(0.325,0.625){\oval(0.24,0.02)[t]}
\put(1.075,0.625){\oval(0.24,0.02)[t]}
\put(0.7,0.625){\oval(0.49,0.02)[t]}
\put(0.7,0.64){\makebox(0,0)[b]{$J_2$}}
\put(0.325,0.64){\makebox(0,0)[b]{$I_2$}}
\put(1.075,0.64){\makebox(0,0)[b]{$I_2$}}
\put(0.45,0.58){\line(0,1){0.04}}
\put(0.95,0.58){\line(0,1){0.04}}
\put(0.45,0.57){\makebox(0,0)[t]{$\frac{1}{8}$}}
\put(0.95,0.57){\makebox(0,0)[t]{$\frac{3}{8}$}}

\put(0.7,0.37){\makebox(0,0)[t]{$\frac{1}{4}$}}
\put(0.2,0.37){\makebox(0,0)[t]{$\frac{1}{8}$}}
\put(1.2,0.37){\makebox(0,0)[t]{$\frac{3}{8}$}}
\put(0.45,0.38){\line(0,1){0.04}}
\put(0.45,0.37){\makebox(0,0)[t]{$\th_0$}}
\put(0.45,0.425){\oval(0.49,0.02)[t]}
\put(0.95,0.425){\oval(0.49,0.02)[t]}
\put(0.45,0.44){\makebox(0,0)[b]{$J_3$}}
\put(0.95,0.44){\makebox(0,0)[b]{$I_3$}}

\put(0.2,0.17){\makebox(0,0)[t]{$\frac{1}{8}$}}
\put(1.2,0.17){\makebox(0,0)[t]{$\frac{1}{4}$}}
\put(0.7,0.17){\makebox(0,0)[t]{$\th_0=\frac{3}{16}$}}
\put(0.325,0.225){\oval(0.24,0.02)[t]}
\put(1.075,0.225){\oval(0.24,0.02)[t]}
\put(0.7,0.225){\oval(0.49,0.02)[t]}
\put(0.7,0.24){\makebox(0,0)[b]{$J_4$}}
\put(0.325,0.24){\makebox(0,0)[b]{$I_4$}}
\put(1.075,0.24){\makebox(0,0)[b]{$I_4$}}
\put(0.45,0.18){\line(0,1){0.04}}
\put(0.95,0.18){\line(0,1){0.04}}
\put(0.45,0.17){\makebox(0,0)[t]{$\frac{5}{32}$}}
\put(0.95,0.17){\makebox(0,0)[t]{$\frac{7}{32}$}}

\put(0.15,0.8){\makebox(0,0)[r]{$k=1$:}}
\put(0.15,0.6){\makebox(0,0)[r]{$k=2$:}}
\put(0.15,0.4){\makebox(0,0)[r]{$k=3$:}}
\put(0.15,0.2){\makebox(0,0)[r]{$k=4$:}}
\end{picture}
\caption{Example of the first four intervals for $\th_0=\frac{3}{16}$.
We have an l-step, a c-step, an l-step and another c-step. All
following steps will be also c-steps.}\label{fig:Intervals}
\end{center}\end{figure}

\begin{Def}
\label{Def:Intervals}
Let $\th_0\in\Theta$ be given.
Start with $J_0=[0,1)$. Let $J_{k-1}=[\th^l_k,\th^r_k)$
and define $d_k=\th^r_k-\th^l_k=2^{-k+1}$.
Then $I_k,J_k\subset J_{k-1}$ are constructed
from $J_{k-1}$ according to the following rules.
\bqa
\label{eq:Interval1}
\th_0\in[\th^l_k,\th^l_k+\mbox{$\frac{3}{8}$}d_k) & \Rightarrow &
J_k=[\th^l_k,\th^l_k+\mbox{$\frac{1}{2}$}d_k),\
I_k=[\th^l_k+\mbox{$\frac{1}{2}$}d_k,\th^r_k),\\
\label{eq:Interval2}
\th_0\in[\th^l_k+\mbox{$\frac{3}{8}$}d_k,\th^l_k+\mbox{$\frac{5}{8}$}d_k)
& \Rightarrow &
J_k=[\th^l_k+\mbox{$\frac{1}{4}$}d_k,\th^l_k+\mbox{$\frac{3}{4}$}d_k),
\\
\nonumber
&&
I_k=[\th^l_k,\th^l_k+\mbox{$\frac{1}{4}$}d_k)\cup
[\th^l_k+\mbox{$\frac{3}{4}$}d_k,\th^r_k),\\
\label{eq:Interval3}
\th_0\in[\th^l_k+\mbox{$\frac{5}{8}$}d_k,\th^r_k)
& \Rightarrow &
J_k=[\th^l_k+\mbox{$\frac{1}{2}$}d_k,\th^r_k),\
I_k=[\th^l_k,\th^l_k+\mbox{$\frac{1}{2}$}d_k).
\eqa
\end{Def}
We call the $k$th step of the interval construction an {\em
l-step} if (\ref{eq:Interval1}) applies, a {\em c-step} if
(\ref{eq:Interval2}) applies, and an {\em r-step} if
(\ref{eq:Interval3}) applies, respectively. Fig.\
\ref{fig:Intervals} shows an example for the interval
construction.

Clearly, this is not the only possible way to define an interval
construction. Maybe the reader wonders why we did not center the
intervals around $\th_0$. In fact, this construction would equally
work for the proof. However, its definition would not be easier,
since one still has to treat the case where
$\th_0$ is located close to the boundary. Moreover, our
construction has the nice property that the interval bounds are
finite binary fractions.

Given the interval construction, we can identify the $\th\in I_k$
with lowest complexity:

\begin{Def}
\label{Def:IntervalOpts}
For $\th_0\in\Theta$ and the interval construction $(I_k,J_k)$, let
\bqan
\th^I_k & = & \arg\min\{\K w(\th):\th\in I_k\cap\Theta\},\\
\th^J_k & = & \arg\min\{\K w(\th):\th\in J_k\cap\Theta\},\und\\
\Delta(k) & = & \max\big\{\K w(\th^I_k)-\K w(\th^J_k),0\big\}.
\eqan
If there is no $\th\in I_k\cap\Theta$, we set $\Delta(k)=\K
w(\th^I_k)=\infty$.
\end{Def}

We can now state the main positive result of this paper. The
detailed proof is deferred to the appendix. Corollaries will be
given in the next section.

\begin{Theorem}
\label{Theorem:UpperBound}
Let $\Theta\subset [0,1]$ be countable,
$\th_0\in\Theta$, and $w_\th=2^{-\K w(\th)}$, where $\K w(\th)$ is
some complexity measure on $\Theta$. Let $\Delta(k)$ be as
introduced in Definition \ref{Def:IntervalOpts} and recall that
$\th^x=\th^{(\alpha,n)}$ depends on $x$'s length and observed
fractions of ones. Then
\beqn
\sum_{n=1}^\infty \Expect (\th_0-\th^x)^2\leqm \K w(\th_0)+
\sum_{k=1}^\infty 2^{-\Delta(k)}\sqrt{\Delta(k)}.
\eeqn
\end{Theorem}

\section{Uniformly Distributed Weights}\label{secUDW}

We are now able to state some positive results following from
Theorem \ref{Theorem:UpperBound}.

\begin{Theorem}
\label{Theorem:Th1}
Let $\Theta\subset[0,1]$ be a countable class of parameters and
$\th_0\in\Theta$ the true parameter. Assume that there are constants
$a\geq 1$ and $b\geq 0$ such that
\beq
\label{eq:Condition1}
\min\big\{\K w(\th):\th\in[\th_0-2^{-k},\th_0+2^{-k}]\cap\Theta,\th\neq\th_0\big\}\geq \frac{k-b}{a}
\eeq
holds for all $k>a\K w(\th_0)+b$. Then we have
\beqn
\sum_{n=1}^\infty \Expect (\th_0-\th^x)^2\leqm a \K w(\th_0)+b\leqm \K w(\th_0).
\eeqn
\end{Theorem}

\begin{Proof}
We have to show that
\beqn
\sum_{k=1}^\infty 2^{-\Delta(k)}\sqrt{\Delta(k)}\leqm a \K w(\th_0)+b,
\eeqn
then the assertion follows from Theorem \ref{Theorem:UpperBound}.
Let $k_1=\lceil a\K w(\th_0)+b+1 \rceil$ and $k'=k-k_1$. Then by
Lemma \ref{Lemma:Iprop} $(iii)$ and (\ref{eq:Condition1}) we have
\bqan
\sum_{k=1}^\infty 2^{-\Delta(k)}\sqrt{\Delta(k)}
& \leq & \sum_{k=1}^{k_1} 1 + \sum_{k=k_1+1}^\infty 2^{-\K
w(\th^I_k)+\K w(\th_0)}
\sqrt{\K w(\th^I_k)-\K w(\th_0)} \\
& \leq & k_1+2^{\K w(\th_0)}\sum_{k=k_1+1}^\infty 2^{-\frac{k-b}{a}}\sqrt{\frac{k-b}{a}}\\
& \leq & k_1+2^{\K w(\th_0)}\sum_{k'=1}^\infty 2^{-\frac{k'+k_1-b}{a}}\sqrt{\frac{k'+k_1-b}{a}}\\
& \leq & a\K w(\th_0)+b+2+\sum_{k'=1}^\infty
2^{-\frac{k'}{a}}\sqrt{\frac{k'}{a}+\K w(\th_0)}.
\eqan
As already seen in the proof of Theorem \ref{Theorem:UpperBound},
$\sqrt{\frac{k'}{a}+\K w(\th_0)}\leq \sqrt{\frac{k'}{a}}+\sqrt{\K w(\th_0)}$,
$\sum_{k'}2^{-\frac{k'}{a}}\leqm a$, and
$\sum_{k'}2^{-\frac{k'}{a}}\sqrt{\frac{k'}{a}}\leqm a$ hold.
The latter is by Lemma \ref{Lemma:IntegralEstimate} $(i)$. This
implies the assertion.
\end{Proof}

Letting $j=\frac{k-b}{a}$, (\ref{eq:Condition1}) asserts that
parameters $\th$ with complexity $\K w(\th)=j$ must have a minimum
distance of $2^{-ja-b}$ from $\th_0$. That is, if parameters with
equal weights are (approximately) uniformly distributed in the
neighborhood of $\th_0$, in the sense that they are not too close
to each other, then fast convergence holds. The next two results
are special cases based on the set of all finite binary fractions,
\beqn
\QQQ_{\BBB^*}=\big\{\th=0.\beta_1\beta_2\ldots \beta_{n-1}1:n\in\NNN,\beta_i\in\BBB\big\}
\cup\big\{0,1\big\}.
\eeqn
If $\th=0.\beta_1\beta_2\ldots\beta_{n-1}1\in\QQQ_{\BBB^*}$, its
length is $l(\th)=n$. Moreover, there is a binary code
$\beta'_1\ldots\beta'_{n'}$ for $n$, having at most
$n'\leq\lfloor\lb (n+1)\rfloor$ bits. Then
$0\beta'_10\beta'_2\ldots 0\beta'_{n'}1\beta_1\ldots \beta_{n-1}$ is a
prefix-code for $\th$. For completeness, we can define the codes
for $\th=0,1$ to be $10$ and $11$, respectively. So we may define
a complexity measure on $\QQQ_{\BBB^*}$ by
\beq
\label{eq:ExampleCoding}
\K w(0)=2,\ \K w(1)=2, \und \K w(\th)=l(\th)+
2\big\lfloor\lb \big(l(\th)+1\big)\big\rfloor \fuer \th\neq 0,1.
\eeq
There are other similar simple prefix codes on $\QQQ_{\BBB^*}$
with the property $\K w(\th)\geq l(\th)$.

\begin{Cor}
\label{Cor:Cor1}
Let $\Theta=\QQQ_{\BBB^*}$, $\th_0\in\Theta$ and $\K w(\th)\geq
l(\th)$ for all $\th\in\Theta$,
and recall $\th^x=\th^{(\alpha,n)}$. Then
$\sum_n \Expect (\th_0-\th^x)^2\leqm \K w(\th_0)$
holds.
\end{Cor}

\begin{Proof}
Condition (\ref{eq:Condition1}) holds with $a=1$ and $b=0$.
\end{Proof}

This is a special case of a uniform distribution of parameters
with equal complexities. The next corollary is more general, it
proves fast convergence if the uniform distribution is distorted
by some function $\ph$.

\begin{Cor}
\label{Cor:Cor2}
Let $\ph:[0,1]\to[0,1]$ be an injective, $N$ times continuously
differentiable function. Let $\Theta=\ph(\QQQ_{\BBB^*})$,
$\K w\big(\ph(t)\big)\geq l(t)$ for all $t\in\QQQ_{\BBB^*}$, and
$\th_0=\ph(t_0)$ for a $t_0\in\QQQ_{\BBB^*}$. Assume that there is
$n\leq N$ and $\eps>0$ such that
\bqan
\left|\frac{d^n\ph}{dt^n}(t)\right|\ \geq \ c \ > & 0 & \for_all t\in[t_0-\eps,t_0+\eps]
\und\\
\frac{d^m\ph}{dt^m}(t_0)  = & 0 & \for_all 1\leq m<n.
\eqan
Then we have
\beqn
\sum \Expect (\th_0-\th^x)^2\leqm n\K w(\th_0)+2\lb(n!)-2\lb c+n\lb\eps \leqm n\K w(\th_0).
\eeqn
\end{Cor}

\begin{Proof}
Fix $j>\K w(\th_0)$, then
\beq
\label{eq:Kw11}
\K w\big(\ph(t)\big)\geq j \for_all
t\in[t_0-2^{-j},t_0+2^{-j}]\cap\QQQ_{\BBB^*}.
\eeq
Moreover, for all
$t\in[t_0-2^{-j},t_0+2^{-j}]$, Taylor's theorem asserts that
\beq
\label{eq:taylor}
\ph(t)=\ph(t_0)+\frac{\frac{d^n\ph}{dt^n}(\tilde t)}{n!}(t-t_0)^n
\eeq
for some $\tilde t$ in $(t_0,t)$ (or $(t,t_0)$ if $t<t_0$). We
request in addition $2^{-j}\leq\eps$, then
$|\frac{d^n\ph}{dt^n}|\geq c$ by assumption. Apply (\ref{eq:taylor}) to
$t=t_0+2^{-j}$ and $t=t_0-2^{-j}$ and define
$k=\lceil jn+\lb(n!)-\lb c\rceil$ in order to obtain
$|\ph(t_0+2^{-j})-\th_0|\geq 2^{-k}$ and
$|\ph(t_0-2^{-j})-\th_0|\geq 2^{-k}$. By injectivity of $\ph$,
we see that $\ph(t)\notin[\th_0-2^{-k},\th_0+2^{-k}]$ if
$t\notin[t_0-2^{-j},t_0+2^{-j}]$. Together with (\ref{eq:Kw11}), this implies
\beqn
\K w(\th)\geq j\geq\frac{k-\lb(n!)+\lb c-1}{n} \for_all
\th\in[\th_0-2^{-k},\th_0+2^{-k}]\cap\Theta.
\eeqn
This is condition (\ref{eq:Condition1}) with $a=n$ and
$b=\lb(n!)-\lb c+1$. Finally, the assumption $2^{-j}\leq\eps$
holds if $k\geq k_1=n\lb\eps+\lb(n!)-\lb c+1$. This gives an
additional contribution to the error of at most $k_1$.
\end{Proof}

Corollary \ref{Cor:Cor2} shows an implication of Theorem
\ref{Theorem:UpperBound} for {\it parameter identification}: A
class of models is given by a set of parameters $\QQQ_{\BBB^*}$
and a mapping $\ph:\QQQ_{\BBB^*}\to\Theta$. The task is to
identify the true parameter $t_0$ or its image $\th_0=\ph(t_0)$.
The injectivity of $\ph$ is not necessary for fast convergence,
but it facilitates the proof. The assumptions of Corollary
\ref{Cor:Cor2} are satisfied if $\ph$ is for example a polynomial.
In fact, it should be possible to prove fast convergence of MDL
for many common parameter identification problems. For sets of
parameters other than $\QQQ_{\BBB^*}$, e.g.\ the set of all
rational numbers
$\QQQ$, similar corollaries can easily be proven.

How large is the constant hidden in ``$\leqm$"? When examining
carefully the proof of Theorem \ref{Theorem:UpperBound}, the
resulting constant is quite huge. This is mainly due to the
frequent ``wasting" of small constants. The sharp bound is
supposably small, perhaps $16$. On the other hand, for the actual
{\it true} expectation (as opposed to its upper bound) and
complexities as in (\ref{eq:ExampleCoding}), numerical simulations
show
$\sum_n \Expect (\th_0-\th^x)^2\leq \frac{1}{2}\K w(\th_0)$.

Finally, we state an implication which almost trivially follows
from Theorem \ref{Theorem:UpperBound} but may be very useful for
practical purposes, e.g. for hypothesis testing (compare
\cite{Rissanen:99}).

\begin{Cor}
\label{Cor:Cor3}
Let $\Theta$ contain $N$ elements, $\K w(\cdot)$ be any complexity
function on $\Theta$, and $\th_0\in\Theta$. Then we have
\beqn
\sum_{n=1}^\infty \Expect (\th_0-\th^x)^2\leqm N+\K w(\th_0).
\eeqn
\end{Cor}

\begin{Proof}
$\sum_k 2^{-\Delta(k)}\sqrt{\Delta(k)}\leq N$ is obvious.
\end{Proof}

\section{The Universal Case} \label{secUC}

We briefly discuss the important universal setup, where $\K
w(\cdot)$ is (up to an additive constant) equal to the prefix
Kolmogorov complexity $\Kpre$ (that is the length of the shortest
self-delimiting program printing $\th$ on some universal Turing
machine). Since $\sum_k 2^{-\Kpre(k)}\sqrt{\Kpre(k)}=\infty$ no
matter how late the sum starts (otherwise there would be a shorter
code for large $k$), Theorem \ref{Theorem:UpperBound} does not
yield a meaningful bound. This means in particular that it does
not even imply our previous result, Theorem \ref{th:Previous}. But
probably the following strengthening of Theorem
\ref{Theorem:UpperBound} holds under the same conditions, which
then easily implies Theorem \ref{th:Previous} up to a constant.

\begin{Conj}
\label{Conj:UpperBound}
$\sum_n \Expect (\th_0-\th^x)^2\leqm \Kpre(\th_0)+
\sum_k 2^{-\Delta(k)}$.
\end{Conj}

Then, take an incompressible finite binary fraction
$\th_0\in\QQQ_{\BBB^*}$, i.e.\
$\Kpre(\th_0)\equa l(\th_0)+\Kpre\big(l(\th_0)\big)$. For $k>l(\th_0)$, we can
reconstruct $\th_0$ and $k$ from $\th^I_k$ and $l(\th_0)$ by just
truncating $\th^I_k$ after $l(\th_0)$ bits. Thus
$\Kpre(\th^I_k)+\Kpre\big(l(\th_0)\big)\geqm
\Kpre(\th_0)+\Kpre\big(k|\th_0,\Kpre(\th_0)\big)$ holds. Using Conjecture
\ref{Conj:UpperBound}, we obtain
\beq\label{polylogbnd}
  \sum_n \Expect (\th_0-\th^x)^2
  \leqm \Kpre(\th_0)+2^{\Kpre(l(\th_0))}
  \leqm l(\th_0)\big(\lb l(\th_0)\big)^2,
\eeq
where the last inequality follows from the example coding given in
(\ref{eq:ExampleCoding}).

So, under Conjecture \ref{Conj:UpperBound}, we obtain a bound
which slightly exceeds the complexity $\Kpre(\th_0)$ if $\th_0$
has a certain structure. It is not obvious if the same holds for
all computable $\th_0$. In order to answer this question positive,
one could try to use something like \cite[Eq.(2.1)]{Gacs:83}. This
statement implies that as soon as $\Kpre(k)\geq K_1$ for all
$k\geq k_1$, we have $\sum_{k\geq k_1} 2^{-\Kpre(k)}\leqm
2^{-K_1}K_1(\lb K_1)^2$. It is possible to prove an analogous
result for $\th^I_k$ instead of $k$, however we have not found an
appropriate coding that does without knowing $\th_0$. Since the
resulting bound is exponential in the code length, we therefore
have not gained anything.

Another problem concerns the size of the multiplicative constant
that is hidden in the upper bound. Unlike in the case of uniformly
distributed weights, it is now of exponential size, i.e.\
$2^{O(1)}$. This is no artifact of the proof, as the following
example shows.

\begin{Example}
\label{Ex:Sensitive} Let $U$ be some universal Turing machine. We
construct a second universal Turing machine $U'$ from U as
follows: Let $N\geq 1$. If the input of $U'$ is $1^Np$, where
$1^N$ is the string consisting of $N$ ones and $p$ is some
program, then $U$ will be executed on $p$. If the input of $U'$ is
$0^N$, then $U'$ outputs $\frac{1}{2}$. Otherwise, if the input of
$U'$ is $x$ with $x\in\BBB^N\setminus\{0^N,1^N\}$, then $U'$
outputs $\frac{1}{2}+2^{-x-1}$. For $\th_0=\frac{1}{2}$, the
conditions of Corollary \ref{Cor:Counterex} are satisfied (where
the complexity is relative to $U'$), thus
$\sum_n\Expect(\th^x-\th_0)^2\geqm 2^N$.
\end{Example}

Can this also happen if the underlying universal Turing machine is
not ``strange" in some sense, like $U'$, but ``natural"? Again
this is not obvious. One would have to define first an appropriate
notion of a ``natural" universal Turing machine which rules out
cases like $U'$. If $N$ is of reasonable size, then one can even
argue that $U'$ {\it is} natural in the sense that its compiler
constant relative to $U$ is small.

There is a relation to the class of all {\it deterministic}
(generally non-i.i.d.) measures. Then MDL predicts the next symbol
just according to the {\it monotone complexity} $\Km$, see
\cite{Hutter:03unimdl}. According to \cite[Theorem
5]{Hutter:03unimdl}, $2^{-\Km}$ is very close to the universal
semimeasure $\M$ \cite{Zvonkin:70,Levin:73random}. Then the
total prediction error (which is defined slightly differently in
this case) can be shown to be bounded by
$2^{O(1)}\Km(x\ltinf)^3$
\cite{Hutter:04unimdlx}. The similarity to the (unproven) bound
(\ref{polylogbnd}) ``huge constant $\times$ polynomial" for the
universal Bernoulli case is evident.

\section{Discussion and Conclusions}\label{secDC}

We have discovered the fact that the instantaneous and the
cumulative loss bounds can be \emph{incompatible}. On the one
hand, the cumulative loss for MDL predictions may be exponential,
i.e. $2^{\K w(\th_0)}$. Thus it implies almost sure convergence at
a slow speed, even for arbitrary discrete model classes
\cite{Poland:04mdl}. On the other hand, the instantaneous loss is
always of order $\frac{1}{n}\K w(\th_0)$, implying fast
convergence in probability and a cumulative loss bound of $\K
w(\th_0)\ln n$. Similar logarithmic loss bounds can be found in
the literature for continuous model classes \cite{Rissanen:96}.

A different approach to assess convergence speed is presented in
\cite{Barron:91}. There, an index of resolvability is introduced,
which can be interpreted as the difference of the expected MDL
code length and the expected code length under the true model. For
discrete model classes, they show that the index of resolvability
converges to zero as $\frac{1}{n}\K w(\th_0)$ \cite[Equation
(6.2)]{Barron:91}. Moreover, they give a convergence of the
predictive distributions in terms of the Hellinger distance
\cite[Theorem 4]{Barron:91}. This implies a cumulative (Hellinger)
loss bound of $\K w(\th_0)\ln n$ and therefore fast convergence in
probability.

If the prior weights are arranged nicely, we have proven a small
finite loss bound $\K w(\th_0)$ for MDL (Theorem
\ref{Theorem:UpperBound}). If parameters of equal complexity are
uniformly distributed or not too strongly distorted (Theorem
\ref{Theorem:Th1} and Corollaries), then the error is within a
small multiplicative constant of the complexity
$\K w(\th_0)$. This may be applied e.g.\ for the case of parameter
identification (Corollary \ref{Cor:Cor2}). A similar result holds
if $\Theta$ is finite and contains only few parameters (Corollary
\ref{Cor:Cor3}), which may be e.g.\ satisfied for hypothesis
testing. In these cases and many others, one can interpret the
conditions for fast convergence as the presence of prior
knowledge. One can show that if a predictor converges to the
correct model, then it performs also well under arbitrarily chosen
bounded loss-functions \cite[Theorem 4]{Hutter:02spupper}.
From an information theoretic viewpoint one may
interpret the conditions for a small bound in Theorem
\ref{Theorem:UpperBound} as ``good codes".

We have proven our positive results only for Bernoulli classes, of
course it would be desirable to cover more general i.i.d.\
classes. At least for \emph{finite} alphabet, our assertions are
likely to generalize, as this is the analog to Theorem
\ref{th:Previous} which also holds for arbitrary finite alphabet.
Proving this seems even more technical than Theorem
\ref{Theorem:UpperBound} and therefore not very interesting. (The
interval construction has to be replaced by a sequence of nested
sets in this case. Compare also the proof of the main result in
\cite{Rissanen:96}.) For small alphabets of size $A$, meaningful
bounds can still be obtained by chaining our bounds $A-1$ times.

It seems more interesting to ask if our results can be
\emph{conditionalized} with respect to inputs. That is, in each
time step, we are given an input and have to predict a label. This
is a standard classification problem, for example a binary
classification in the Bernoulli case. While it is straightforward
to show that Theorem \ref{th:Previous} still holds in this setup
\cite{Poland:05mdlreg}, it is not clear in which way the present
proofs can be adapted. We leave this interesting question open.

We conclude with another open question. In abstract terms, we have
proven a convergence result for the Bernoulli case by mainly
exploiting the {\it geometry} of the space of distributions. This
has been quite easy in principle, since for Bernoulli this space
is just the unit interval (for i.i.d it is the space of
probability vectors). It is not at all obvious if this approach
can be transferred to general (computable) measures.

\begin{appendix}
\section{Proof of Theorem \ref{Theorem:UpperBound}}

The proof of Theorem  \ref{Theorem:UpperBound} requires some
preparations. We start by showing assertions on the interval
construction from Definition \ref{Def:Intervals}.

\begin{Lemma}
\label{Lemma:Iprop}
The interval construction has the following properties.
\bqan
(i) && |J_k| = 2^{-k},\\
(ii) && d(\th_0,I_k) \geq 2^{-k-2},\\
(iii) && \max_{\th\in I_k}|\th-\th_0| \leq 2^{-k+1},\\
(iv) && d(J_{k+5},I_k) \geq 15\cdot 2^{-k-6}.
\eqan
\end{Lemma}

By $d(\cdot,\cdot)$ we mean the Euclidean distance:
$d(\tilde\th,I)=\min\{|\tilde\th-\th|:\th\in I\}$ and
$d(J,I)=\min\{d(\tilde\th,I):\tilde\th\in J\}$.

\begin{Proof}
The first three equations are fairly obvious. The last estimate
can be justified as follows. Assume that $k$th step of the
interval construction is a c-step, the same argument applies if it
is an l-step or an r-step. Let $c$ be the center of $J_k$ and
assume without loss of generality
$\th_0\leq c$. Define $\th_I=\max\{\th\in I_k:\th<c\}$ and
$\th_J=\min\{\th\in J_{k+5}\}$ (recall the general assumption
$\th\in\Theta$ for all $\th$ that occur, i.e.
$\th_I,\th_J\in\Theta$). Then $\th_I=c-2^{-k-1}$ and
$\th_J\geq c-2^{-k-2}-2^{-k-6}$, where equality holds if $\th_0=c-2^{-k-2}$.
Consequently, $\th_J-\th_I\geq 2^{-k-1}-2^{-k-2}-2^{-k-6}=15\cdot 2^{-k-6}$.
This establishes the claim.
\end{Proof}

Next we turn to the minimum complexity elements in the intervals.

\begin{Prop}
\label{Prop:IntervalOpts}
The following assertions hold for all $k\geq 1$.
\bqan
(i) && \K w(\th^J_k) \leq \K w(\th_0), \\
(ii) && \K w(\th^J_{k+6}) \geq \K w(\th^J_{k}), \\
(iii) && \K w(\th^I_{k+1}) \geq \K w(\th^J_{k}), \\
(iv) && \sum_{k=1}^\infty \max\big\{\K w(\th^J_{k+5})-\K w(\th^I_k),0\big\} \leq 6 \K w(\th_0), \\
\eqan
\end{Prop}

\begin{Proof}
The first three inequalities follow from $\th_0\in J_{k}$ and
$J_{k+6},I_{k+1}\subset J_{k}$. This implies
\bqan
\lefteqn{\sum_{j=0}^m \max\big\{\K w(\th^J_{6j+6})-\K w(\th^I_{6j+1}),0\big\}}\\
& \leq &
\max\big\{\K w(\th^J_{6})-\K w(\th^I_{1}),0\big\} + \sum_{j=1}^m
\max\big\{\K w(\th^J_{6j+6})-\K w(\th^J_{6j}),0\big\}\\
& \leq & \K w(\th^J_{6}) + \sum_{j=1}^m\big[
\K w(\th^J_{6j+6})-\K w(\th^J_{6j})\big] = \K w(\th^J_{6m+6}) \leq \K w(\th_0)
\eqan
for all $m\geq 0$. By the same argument, we have
\beqn
\sum_{j=0}^m \max\big\{\K w(\th^J_{6j+i+5})-\K w(\th^I_{6j+i}),0\big\}\leq \K w(\th_0)
\eeqn
for all $1\leq i\leq 6$ (use $(iii)$ in the first inequality,
$(ii)$ in the second, and $(i)$ in the last). This implies $(iv)$.
Clearly, we could everywhere substitute $5$ by some constant $k'$
and $6$ by $k'+1$, but we will need the assertion only for the
special case.
\end{Proof}

Consider the case that $\th_0$ is located close to the boundary of
$[0,1]$. Then the interval construction involves for long time only
l-steps, if we assume without loss of generality
$\th_0\leq\frac{1}{2}$. We will need to treat this case
separately, since the estimates for the general situation work
only as soon as at least one c-step has taken place. Precisely,
the interval construction consists only of l-steps as long as
\beqn
\mbox{$
\th_0<\frac{3}{4}2^{-k}$, i.e.\ $k<-\lb\th_0+\lb(\frac{3}{4})$.}
\eeqn
We therefore define
\beq
\label{eq:k0}
\mbox{$
k_0=\max\big\{0,\lfloor-\lb\th_0+\lb\frac{3}{4}\rfloor\big\}
$}
\eeq
and observe that the $(k_0+1)$st step is the first c-step.
We are now prepared to give the main proof.

\textbf{Proof} of Theorem \ref{Theorem:UpperBound}.
Assume $\th_0\in\Theta\setminus\{0,1\}$, the case
$\th_0\in\{0,1\}$ is handled like Case 1a below
and will be left to the reader.

Before we start, we will show that the contribution of $\th=1$ to
the total error is bounded by $\frac{1}{4}$. This is immediate, since $1$
cannot become the maximizing element as soon as $x\neq 1^n$.
Therefore the contribution is bounded by
\beq
\label{eq:ContributionBoundary}
\sum_{n=1}^\infty (1-\th_0)^2 p(1^n)=
(1-\th_0)^2 \sum_{n=1}^\infty  \th_0^n = \th_0(1-\th_0)\leq\mbox{$\frac{1}{4}$}.
\eeq
The same is true for the contribution of $\th=0$.

As already mentioned, we first estimate the contributions of
$\th\in I_k$ for small $k$ if the true parameter $\th_0$ is
located close to the boundary. To this aim, we assume
$\th_0\leq\frac{1}{2}$ without loss of generality. We know that
the interval construction involves only l-steps as long as $k\leq
k_0$, see (\ref{eq:k0}). The very last five of these $k$ still
require a particular treatment, so we start with $k\leq k_0-5$ and
$\alpha$ is far from $\th_0$. (If $k_0-5<1$, then there is nothing
to estimate.)

{\bf Case 1a:} $k\leq k_0-5$, $j\leq k_1$, $\alpha\in
I_j=[2^{-j},2^{-j+1})$, where $k_1=k+\lceil\lb (k_0-k-3)\rceil+2$.
The probability of $\alpha$ does not exceed $p(2^{-j})$. The
squared error may clearly be upper bounded by
$2^{-2k+2}=O(2^{-2k})$. For $n<2^j$, no such fractions can occur,
so we may consider only $n=2^j+n'$, $n'\geq 0$. Finally, there are
at most $\lceil n\cdot 2^{-j-1}\rceil=O(2^{-j}n)$ fractions
$\alpha\in I_j$. This follows from the general fact that if
$I\subset(0,1)$ is any half-open or open interval of length at
most $l$, then at most $\lceil n l\rceil$ observed fractions can
be located in $I$.

We now derive an estimate for the probability which is
\beqn
p(\alpha|n)\leq p(2^{-j}|n)\leqm n^{-\frac{1}{2}}
2^{\frac{j}{2}}\exp\big[-n\cdot D(2^{-j}\|\th_0)\big]
\eeqn
according to Lemma \ref{Lemma:BinomialBounds}. Then, Lemma
\ref{Lemma:EntropyIneq} (v) implies
\beqn
\exp\big[-n D(2^{-j}\|\th_0)\big]
\leq\exp\big[-(2^j+n') D(2^{-j}\|2^{-k_0})\big]
\leq \exp\big[n'2^{-j}(k_0-j-1)\big].
\eeqn
Taking into account the upper bound for the squared error $O(2^{-2k})$ and
the maximum number of fractions $O(2^{-j}n)$,
the contribution $C(k,j)$ can be upper bounded by
\beqn
C(k,j)\leqm
\sum_{n=2^j}^\infty p(2^{-j}|n)2^{-2k}\cdot 2^{-j}n
\leqm
\sum_{n'=0}^\infty
2^{-2k-\frac{j}{2}}\sqrt n\cdot
\exp\big[n'2^{-j}(k_0-j-1)\big].
\eeqn
Decompose the right hand side using $\sqrt{n}\leq\sqrt{2^j}+\sqrt{n'}$.
Then we have
\bqan
\sum_{n'=0}^\infty
2^{-2k-\frac{j}{2}}\sqrt{2^j}\cdot \exp\big[n'2^{-j}(k_0-j-1)\big]
& \leqm & 2^{-2k+j}(k_0-j-1)^{-1} \und\\
\sum_{n'=0}^\infty
2^{-2k-\frac{j}{2}}\sqrt{n'}\cdot \exp\big[n'2^{-j}(k_0-j-1)\big]
& \leqm & 2^{-2k+j}(k_0-j-1)^{-\frac{3}{2}}
\eqan
where the first inequality is straightforward and the second holds
by Lemma \ref{Lemma:IntegralEstimate} $(i)$. Letting $k'=k_0-k-3$,
we have $k'\geq 2$ and
\beqn
(k_0-j-1)^{-\frac{3}{2}} \leq (k_0-j-1)^{-1}
\leq (k_0-k_1-1)^{-1}
= \big(k'-\lceil\lb k'\rceil\big)^{-1}.
\eeqn
Thus we may conclude
\bqa
\label{eq:Ckleq1}
C(k,\,\leq\!k_1) & := & \sum_{j=1}^{k_1}C(k,j)\ \leqm \
\sum_{j=1}^{k+\lceil\lb k'\rceil+2}\frac{2^{-2k+j}}{k'-\lceil\lb k'\rceil}\\
\nonumber
& \leqm & 2^{-k}\frac{k'}{k'-\lceil\lb k'\rceil}
\ \leq\
2^{-k}\left(1+\frac{\lceil\lb k'\rceil}{k'-\lceil\lb k'\rceil}\right)
\ \leq\ 3\cdot 2^{-k}
\eqa
(the last inequality is sharp for $k'=3$).

{\bf Case 1b:} $k\leq k_0-5$, $\alpha\leq 2^{-k_1}$
(recall $k_1=k+\lceil\lb (k_0-k-3)\rceil+2$).
This means that we consider $\alpha$ close to $\th_0$.
By (\ref{eq:Beating})
we know that $\th_0$ beats $\th\in I_k$ if
\beqn
n\cdot D^\alpha\big(\th_0\|\th\big)\geq \ln 2\big(\K w(\th_0)-\K
w(\th)\big)
\eeqn
holds. This happens certainly for $n\geq N_1:=\ln 2 \cdot \K
w(\th_0)\cdot 2^{k+4}$, since Lemma \ref{Lemma:ForTheorem3} below
asserts
$D^\alpha\big(\th_0\|\th\big)\geq 2^{-2-4}$.
Thus only smaller $n$ can contribute.
The total probability of all $\alpha\leq 2^{-k_1}$ is clearly
bounded by means of
\beqn
\sum_{\alpha} p(\alpha|n)\leq 1.
\eeqn
The jump size, i.e.\ the squared error, is again $O(2^{-2k})$. Hence
the total contribution caused in $I_k$ by $\alpha\leq 2^{-k_1}$ can
thus be upper bounded by
\beqn
C(k,\,>\!k_1)\leqm \sum_{n=1}^{N_1} 2^{-2k}\leqm \K
w(\th_0)2^{-k},
\eeqn
where $C(k,\,>\!k_1)$ is the obvious abbreviation for this
contribution. Together with (\ref{eq:Ckleq1}) this implies
$C(k)\leqm \K w(\th_0) 2^{-k}$ and therefore
\beq
\sum_{k=1}^{k_0-5} C(k)\leqm \K w(\th_0).
\label{eq:Estimate1}
\eeq

This finishes the estimates for $k\leq k_0-5$.
We now will consider the indices
\beqn
k_0-4\leq k\leq k_0
\eeqn
and show that the contributions caused by these $\th\in I_k$ is at
most $O\big(\K w(\th_0)\big)$.

{\bf Case 2a:} $k_0-4\leq k\leq k_0$, $j\leq k+5$, $\alpha\in I_j$.
Assume that $\th\in I_k$ starts contributing only for $n>n_0$.
This is not relevant here, and we will set $n_0=0$ for the moment,
but then we can reuse the following computations later.
Consequently we have $n=n_0+n'$, and from Lemma \ref{Lemma:BinomialBounds}
we obtain
\beq
\label{eq:C1b1}
p(\alpha|n) \leqm n^{-\frac{1}{2}}
2^{\frac{k_0}{2}}\exp\big[-(n_0+n')\cdot D(\alpha\|\th_0)\big].
\eeq
Lemma \ref{Lemma:Iprop} implies $d(\alpha,\th_0)\geq 2^{-j-2}$ and thus
\beq
\label{eq:C1b11}
D(\alpha\|\th_0)\geq \frac{2^{-2j-4}}{2\cdot 2^{-k_0}}
= 2^{-2j-5+k_0}.
\eeq
according to Lemma \ref{Lemma:EntropyIneq} (iii).
Therefore we obtain
\beq
\label{eq:C1b2}
\exp\big[-(n_0+n')\cdot D(\alpha\|\th_0)\big]\leq
\exp\big[-n_0\cdot D(\alpha\|\th_0)\big]
\exp\big[-n'2^{-2j-5+k_0}\big].
\eeq
Again the maximum square error is $O(2^{-2k})$, the maximum number
of fractions is $O(n2^{-j})$. Therefore
\beq
\label{eq:C1b3}
C(k,j)\leqm
\exp\big[-n_0D(\alpha\|\th_0)\big]
\sum_{n'=1}^\infty 2^{-2k-j+\frac{k_0}{2}}\sqrt{n_0+n'}
\exp\big[-n'2^{-2j-5+k_0}\big].
\eeq
We have
\bqa
\label{eq:C1b4}
\sum_{n'=1}^\infty 2^{-2k-j+\frac{k_0}{2}}
\exp\big[-n'2^{-2j-5+k_0}\big] & \leqm & 2^{-2k+j-\frac{k_0}{2}}
\leq 2^{-2k+j} \und\quad\\
\label{eq:C1b5}
\sum_{n'=1}^\infty 2^{-2k-j+\frac{k_0}{2}}\sqrt{n'}
\exp\big[-n'2^{-2j-5+k_0}\big] & \leqm & 2^{-2k+2j-k_0}
\leq 2^{-2k+2j},
\eqa
where the first inequality is straightforward and the second
follows from Lemma \ref{Lemma:IntegralEstimate} $(i)$. Observe
$\sum_{j=1}^{k+5}2^{j}\leqm 2^{k}$, $\sum_{j=1}^{k+5}2^{2j}\leqm
2^{2k}$, and $\sqrt{n}\leq\sqrt{n_0}+\sqrt{n'}$ in order to obtain
\beq
\label{eq:C1b6}
C(k,\,\leq\!k+5)\leqm
\exp\big[-n_0D(\alpha\|\th_0)\big]
(1+2^{-k}\sqrt{n_0}).
\eeq
The right hand side depends not only on $k$ and $n_0$, but formally also on
$\alpha$ and even on $\th$, since $n_0$ itself depends on $\alpha$
and $\th$.
Recall that for this case we have agreed on $n_0=0$, thus $C(k,\,\leq\!k+5)=O(1)$.

{\bf Case 2b:} $k_0-4\leq k\leq k_0$, $\alpha\in J_{k+5}$.
As before, we will argue that then
$\th\in I_k$ can be the maximizing element only for small $n$.
Namely, $\th_0$ beats $\th$ if
$n\cdot D^\alpha\big(\th_0\|\th\big)\geq \ln 2\big(\K w(\th_0)-\K w(\th)\big)$
holds.
Since $D^\alpha\big(\th_0\|\th\big)\geq 2^{-2k-5}$ as stated in Lemma
\ref{Lemma:ForTheorem3} below, this happens certainly for
$n\geq N_1:=\ln 2 \cdot \K w(\th_0)\cdot 2^{2k+5}$,
thus only smaller $n$ can contribute. Note that in order to
apply Lemma \ref{Lemma:ForTheorem3}, we need $k\geq k_0-4$.
Again the total probability of all $\alpha$ is at most $1$ and the
jump size is $O(2^{-2k})$, hence
\beqn
C(k,\,>\!k+5)\leqm\sum_{n=1}^{N_1} 2^{-2k}\leqm \K w(\th_0).
\eeqn
Together with $C(k,\,\leq\!k+5)=O(1)$ this implies $C(k)\leqm \K
w(\th_0)$ and thus
\beq
\sum_{k=k_0-4}^{k_0} C(k)\leqm \K w(\th_0).
\label{eq:Estimate2}
\eeq

This completes the estimate for the initial l-steps.
We now proceed with the main part of the proof. At this point, we drop
the general assumption $\th_0\leq\frac{1}{2}$, so that
we can exploit the symmetry otherwise if convenient.

{\bf Case 3a:} $k\geq k_0+1$, $j\leq k+5$, $\alpha\in I_j$.
For this case, we may repeat the computations
(\ref{eq:C1b1})-(\ref{eq:C1b6}), arriving at
\beq
\label{eq:C1main}
C(k,\,\leq\!k+5)\leqm
\exp\big[-n_0D(\alpha\|\th_0)\big]
(1+2^{-k}\sqrt{n_0}).
\eeq
The right hand side of (\ref{eq:C1main}) depends on $k$ and $n_0$ and formally
also on $\alpha$ and $\th$.
We now come to the crucial point of this proof:
\beqn
\mbox{\em{For most $k$, $n_0$ is considerably larger than $0$.}}
\eeqn
That is, for most $k$, $\th\in I_k$ starts contributing
late, i.e.\ for large $n$. This will cause the right hand side
of (\ref{eq:C1main}) to be small.

We know that $\th_0$ beats $\th\in I_k$
for {\em any} $\alpha\in[0,1]$ as long as
\beq
\label{eq:E2}
n D^\alpha(\th\|\th_0)\leq \ln 2\big(\K w(\th)-\K w(\th_0)\big)
\eeq
holds. We are interested in for which $n$ this must happen
regardless of $\alpha$, so assume that $\alpha$ is close enough to $\th$
to make $D^\alpha(\th\|\th_0)>0$. Since
$\K w(\th)\geq \K w(\th^I_k)$, we see that
(\ref{eq:E2}) holds if
\beqn
n\leq n_0(k,\alpha,\th):=
\frac{\ln 2\cdot \Delta(k)}{D^\alpha(\th\|\th_0)}.
\eeqn
We show the following two relations:
\bqa
\label{eq:toshow1}
\exp\big[-n_0(k,\alpha,\th)D(\alpha\|\th_0)\big] & \leq &
2^{-\Delta(k)} \und\\
\label{eq:toshow2}
\exp\big[-n_0(k,\alpha,\th)D(\alpha\|\th_0)\big]2^{-k}\sqrt{n_0(k,\alpha,\th)} & \leqm &
2^{-\Delta(k)}\sqrt{\Delta(k)},
\eqa
regardless of $\alpha$ and $\th$. Since
$D(\alpha\|\th_0)\geq D(\alpha\|\th_0)-D(\alpha\|\th)=D^\alpha(\th\|\th_0)$,
(\ref{eq:toshow1}) is immediate.
In order to verify (\ref{eq:toshow2}), we observe that
\beqn
D(\alpha\|\th_0)\geq 2^{-2j-5+k_0}\geq 2^{-2k-15+k_0} \geq 2^{-2k-15}
\eeqn
holds as in (\ref{eq:C1b11}).
So for those $\alpha$ and $\th$ having
\beq
\label{eq:alpha1}
\eta:=\frac{2^{-2k-15}}{D^\alpha(\th\|\th^J_{k+5})}\geq 1,
\eeq
we obtain
\bqan
\exp\big[-n_0(k,\alpha,\th)D(\alpha\|\th_0)\big]2^{-k}\sqrt{n_0(k,\alpha,\th)} & \leq &
2^{-\Delta(k)\eta}
2^{-k}\sqrt{\ln 2\cdot \Delta(k)\eta 2^{2k+15}}
\\ & \leqm & 2^{-\Delta(k)}\sqrt{\Delta(k)}.
\eqan
since $\eta\geq 1$.
If on the other hand (\ref{eq:alpha1}) is not
valid, then $D^\alpha(\th\|\th^J_{k+5})\leqm 2^{-2k}$ holds, which
together with $D(\alpha\|\th_0)\geq D^\alpha(\th\|\th_0)$
again
implies (\ref{eq:toshow2}).

So we conclude that the dependence on $\alpha$ and $\th$ of the right
hand side of (\ref{eq:C1main}) is indeed only a formal one. So we
obtain
$C(k,\,\leq\!k+5)\leqm 2^{-\Delta(k)}\sqrt{\Delta(k)}$,
hence
\beq
\label{eq:Estimate3a}
\sum_{k=k_0+1}^\infty C(k,\,\leq\!k+5)\leqm
\sum_{k=1}^\infty 2^{-\Delta(k)}\sqrt{\Delta(k)}.
\eeq

{\bf Case 3b:} $k\geq k_0+1$, $\alpha\in J_{k+5}$.
We know that $\th^J_{k+5}$ beats $\th$ if
\beqn
n\geq \ln 2 \cdot \max\big\{\K w(\th^J_{k+5})-\K
w(\th),0\big\}\cdot 2^{2k+5},
\eeqn
since $D^\alpha\big(\th^J_{k+5}\|\th\big)\geq 2^{-2k-5}$ according
to Lemma \ref{Lemma:ForTheorem3}. Since $\K w(\th)\geq
\K w(\th^I_k)$, this happens certainly for $n\geq N_1:=\ln 2 \cdot
\max\big\{\K w(\th^J_{k+5})-\K w(\th^I_k),0\big\}\cdot 2^{2k+5}$. Again
the total probability of all $\alpha$ is at most $1$ and the jump
size is $O(2^{-2k})$. Therefore we have
\beqn
C(k,\,>\!k+5)\leqm\sum_{n=1}^{N_1} 2^{-2k}\leqm
\max\big\{\K w(\th^J_{k+5})-\K w(\th^I_k),0\big\}.
\eeqn
Using Proposition \ref{Prop:IntervalOpts} $(iv)$, we conclude
\beq
\label{eq:Estimate3b}
\sum_{k=k_0+1}^\infty C(k,\,>\!k+5)\leqm \K w(\th_0).
\eeq

Combining all estimates for $C(k)$, namely
(\ref{eq:Estimate1}), (\ref{eq:Estimate2}),
(\ref{eq:Estimate3a}) and (\ref{eq:Estimate3b}), the assertion follows. \hfill $\Box$

\begin{Lemma}
\label{Lemma:ForTheorem1}
Let $1\leq k\leq k_0-5$, $k_1=k+\lceil\lb (k_0-k-3)\rceil+2$,
$\th\geq 2^{-k}$, and $\alpha\leq 2^{-k_1}$. Then
$D^\alpha(\th_0\|\th)\geq 2^{-k-4}$ holds.
\end{Lemma}

\begin{Proof}
By Lemma \ref{Lemma:EntropyIneq} $(iii)$ and $(vii)$, we have
\bqan
D(\alpha\|\th) & \geq & D(2^{-k_1}\|2^{-k})
\ \geq \ \frac{\big(2^{-k}-2^{-k_1}\big)^2}{2\cdot 2^k(1-2^k)}\\
& \geq & 2^{-k-1}\big(1-2^{-\lceil\lb (k_0-k-3)\rceil-2}\big)
\ \geq \ 7\cdot 2^{-k-4} \und\\[2ex]
D(\alpha\|\th_0) & \leq & D(2^{-k_1}\|2^{-k_0-1})
\ \leq \ 2^{-k_1}(k_0+1-k_1)\\
& \leq & 2^{-k-2}\frac{k_0-k-\lceil\lb (k_0-k-3)\rceil-1}{k_0-k-3}
\ \leq \ 6\cdot 2^{-k-4}
\eqan
(the last inequality is sharp for $k=k_0-5$).
This implies $D^\alpha(\th_0\|\th)=D(\alpha\|\th)-D(\alpha\|\th_0)
\geq 2^{-k-4}$.
\end{Proof}

\begin{Lemma}
\label{Lemma:ForTheorem3}
Let $k\geq k_0-4$,
$\th\in I_k$, and $\alpha,\tilde\th\in J_{k+5}$. Then we have
$D^\alpha(\tilde\th\|\th)\geq 2^{-2k-5}$.
\end{Lemma}

\begin{Proof}
Assume $\th\leq\frac{1}{2}$ without loss of generality. Moreover,
we will only present the case $\tilde\th\leq\th\leq\frac{1}{4}$,
the other cases are similar and simpler. From
Lemma \ref{Lemma:EntropyIneq} $(iii)$ and $(iv)$ and
Lemma \ref{Lemma:Iprop} we know that
\bqan
D(\alpha\|\th) & \geq & \frac{(\alpha-\th)^2}{2\th(1-\th)}
\ \geq \ \frac{15^2 2^{-2k-12}}{2\th} \und\\
D(\alpha\|\tilde\th) & \leq & \frac{3(\alpha-\tilde\th)^2}{2\alpha(1-\alpha)}
\ \leq \ \frac{4\cdot 3\cdot 2^{-2k-14}}{3\cdot 2\alpha}
\ \leq \ \frac{2\cdot 128\cdot 2^{-2k-14}}{\th}.
\eqan
Note that in order to apply Lemma \ref{Lemma:EntropyIneq} $(iv)$
in the second line we need to know that for $k+5$ a c-step
has already taken place, and the last estimate follows from
$\th\leq 128\alpha$ which is a consequence of $k\geq k_0-4$.
Now the assertion follows from
$D^\alpha(\tilde\th\|\th)=D(\alpha\|\th)-D(\alpha\|\tilde\th)
\geq 2^{-2k-6}\big(15^2 2^{-7}-1\big)\th^{-1}\geq 2^{-2k-5}$.
\end{Proof}

\end{appendix}



\begin{thebibliography}{Hut03b}

\bibitem[BC91]{Barron:91}
A.~R. Barron and T.~M. Cover.
\newblock Minimum complexity density estimation.
\newblock {\em IEEE Trans. on Information Theory}, 37(4):1034--1054, 1991.

\bibitem[BRY98]{Barron:98}
A.~R. Barron, J.~J. Rissanen, and B.~Yu.
\newblock The minimum description length principle in coding and modeling.
\newblock {\em IEEE Trans. on Information Theory}, 44(6):2743--2760, 1998.

\bibitem[CB90]{Clarke:90}
B.~S. Clarke and A.~R. Barron.
\newblock Information-theoretic asymptotics of {Bayes} methods.
\newblock {\em IEEE Trans. on Information Theory}, 36:453--471, 1990.

\bibitem[G{\'a}c83]{Gacs:83}
P.~G{\'a}cs.
\newblock On the relation between descriptional complexity and algorithmic
  probability.
\newblock {\em Theoretical Computer Science}, 22:71--93, 1983.

\bibitem[GL04]{Gruenwald:04}
P.~Gr{\"u}nwald and J.~Langford.
\newblock Suboptimal behaviour of {B}ayes and {MDL} in classification under
  misspecification.
\newblock In {\em 17th Annual Conference on Learning Theory (COLT)}, pages
  331--347, 2004.

\bibitem[Hut01]{Hutter:01alpha}
M.~Hutter.
\newblock Convergence and error bounds for universal prediction of nonbinary
  sequences.
\newblock {\em Proc. 12th Eurpean Conference on Machine Learning (ECML-2001)},
  pages 239--250, December 2001.

\bibitem[Hut03a]{Hutter:02spupper}
M.~Hutter.
\newblock Convergence and loss bounds for {Bayesian} sequence prediction.
\newblock {\em IEEE Trans. on Information Theory}, 49(8):2061--2067, 2003.

\bibitem[Hut03b]{Hutter:03optisp}
M.~Hutter.
\newblock Optimality of universal {B}ayesian prediction for general loss and
  alphabet.
\newblock {\em Journal of Machine Learning Research}, 4:971--1000, 2003.

\bibitem[Hut03c]{Hutter:03unimdl}
M.~Hutter.
\newblock Sequence prediction based on monotone complexity.
\newblock In {\em Proc. 16th Annual Conference on Learning Theory
  ({COLT-2003})}, Lecture Notes in Artificial Intelligence, pages 506--521,
  Berlin, 2003. Springer.

\bibitem[Hut04]{Hutter:04unimdlx}
M.~Hutter.
\newblock Sequential predictions based on algorithmic complexity.
\newblock {\em Journal of Computer and System Sciences}, 2005.
\newblock 72(1):95--117.

\bibitem[Lev73]{Levin:73random}
L.~A. Levin.
\newblock On the notion of a random sequence.
\newblock {\em Soviet Math. Dokl.}, 14(5):1413--1416, 1973.

\bibitem[Li99]{Li:99}
J.~Q. Li.
\newblock {\em Estimation of Mixture Models}.
\newblock PhD thesis, Dept. of Statistics. Yale University, 1999.

\bibitem[LV97]{Li:97}
M.~Li and P.~M.~B. Vit\'anyi.
\newblock {\em An introduction to {Kolmogorov} complexity and its
  applications}.
\newblock Springer, 2nd edition, 1997.

\bibitem[PH04a]{Poland:04mdl}
J.~Poland and M.~Hutter.
\newblock Convergence of discrete {MDL} for sequential prediction.
\newblock In {\em 17th Annual Conference on Learning Theory (COLT)}, pages
  300--314, 2004.

\bibitem[PH04b]{Poland:04mdlspeed}
J.~Poland and M.~Hutter.
\newblock On the convergence speed of {MDL} predictions for {B}ernoulli
  sequences.
\newblock In {\em International Conference on Algorithmic Learning Theory
  (ALT)}, pages 294--308, 2004.

\bibitem[PH05]{Poland:05mdlreg}
J.~Poland and M.~Hutter.
\newblock Strong asymptotic assertions for discrete {MDL} in regression and
  classification.
\newblock In {\em Benelearn 2005 (Ann. Machine Learning Conf. of Belgium and
  the Netherlands)}, 2005.

\bibitem[Ris96]{Rissanen:96}
J.~J. Rissanen.
\newblock Fisher {I}nformation and {S}tochastic {C}omplexity.
\newblock {\em IEEE Trans. on Information Theory}, 42(1):40--47, January 1996.

\bibitem[Ris99]{Rissanen:99}
J.~J. Rissanen.
\newblock Hypothesis selection and testing by the {MDL} principle.
\newblock {\em The Computer Journal}, 42(4):260--269, 1999.

\bibitem[Sol78]{Solomonoff:78}
R.~J. Solomonoff.
\newblock Complexity-based induction systems: comparisons and convergence
  theorems.
\newblock {\em IEEE Trans. Information Theory}, IT-24:422--432, 1978.

\bibitem[VL00]{Vitanyi:00}
P.~M. Vit{\'a}nyi and M.~Li.
\newblock Minimum description length induction, {B}ayesianism, and {K}olmogorov
  complexity.
\newblock {\em IEEE Trans. on Information Theory}, 46(2):446--464, 2000.

\bibitem[Vov97]{Vovk:97}
V.~G. Vovk.
\newblock Learning about the parameter of the {B}ernoulli model.
\newblock {\em Journal of Computer and System Sciences}, 55:96--104, 1997.

\bibitem[Zha04]{Zhang:04}
T.~Zhang.
\newblock On the convergence of {MDL} density estimation.
\newblock In {\em Proc. 17th Annual Conference on Learning Theory (COLT)},
  pages 315--330, 2004.

\bibitem[ZL70]{Zvonkin:70}
A.~K. Zvonkin and L.~A. Levin.
\newblock The complexity of finite objects and the development of the concepts
  of information and randomness by means of the theory of algorithms.
\newblock {\em Russian Mathematical Surveys}, 25(6):83--124, 1970.

\end{thebibliography}
\end{document}